\documentclass[noams,lp]{compositio}
\usepackage{amsbsy,amssymb,amscd,amsfonts,latexsym,amstext,delarray, amsmath,stackrel,lineno,tabularx,url,MnSymbol,tikz,cite}

\input xypic
\usepackage[pdfborder={0 0 0}]{hyperref}

\hypersetup{
    colorlinks = true,
    linkcolor = blue,
    anchorcolor = red,
    citecolor = green,
    filecolor = red,
    pagecolor = red,
    urlcolor = magenta
}
\newcommand*\patchAmsMathEnvironmentForLineno[1]{%
  \expandafter\let\csname old#1\expandafter\endcsname\csname #1\endcsname
  \expandafter\let\csname oldend#1\expandafter\endcsname\csname end#1\endcsname
  \renewenvironment{#1}%
     {\linenomath\csname old#1\endcsname}%
     {\csname oldend#1\endcsname\endlinenomath}}%
\newcommand*\patchBothAmsMathEnvironmentsForLineno[1]{%
  \patchAmsMathEnvironmentForLineno{#1}%
  \patchAmsMathEnvironmentForLineno{#1*}}%
\AtBeginDocument{%
\patchBothAmsMathEnvironmentsForLineno{equation}%
\patchBothAmsMathEnvironmentsForLineno{align}%
\patchBothAmsMathEnvironmentsForLineno{flalign}%
\patchBothAmsMathEnvironmentsForLineno{alignat}%
\patchBothAmsMathEnvironmentsForLineno{gather}%
\patchBothAmsMathEnvironmentsForLineno{multline}%
}
\newtheorem{thm}{Theorem}[section] % number like 3.1, 3.2, 3.3, etc.
\newtheorem*{theorem}{Theorem}
\newtheorem{defn}[thm]{Definition} % numbered with thm
\newtheorem{prop}[thm]{Proposition}
\newtheorem{cor}[thm]{Corollary}
\newtheorem{lem}[thm]{Lemma}

\newtheorem{rem}[thm]{Remark}

\def\Aut{{\rm Aut}}

\def\End{{\rm End}}

\def\Hom{{\rm Hom}}
\def\id{{\rm id}}

\def\mmod{{\rm Mod}}

\def\F{{\mathbb F}}

\def\N{{\mathbb N}}

\def\Q{{\mathbb Q}}
\def\R{{\mathbb R}}
\def\Z{{\mathbb Z}}
\def\H{{\mathbb H}}

\def\cC{{\mathcal C}}

\def\cG{{\mathcal G}}

\def\cF{{\mathcal F}}

\def\cP{{\mathcal P}}

\def\qqq{\,,\,~\forall}

\newcommand{\ie}{{\it i.e.\/}\ }
\newcommand{\eg}{{\it e.g.\/}\ }
\newcommand{\cf}{{\it cf.}}
\newcommand{\opcit}{{\it op.cit.\/}\ }

\def\mod{semimodule}
\def\Mod{{\rm Mod}}
\def\Ab{{\mathfrak{Ab}}}

\def\mods{semimodules}

\def\id{{\mbox{Id}}}

\def\Hom {{\mbox{Hom}}}

\def\End{{\mbox{End}}}

\def\fr{{\rm Fr}}

\def\Se{\frak{ Sets}}
\def\fin{\frak{ Fin}}

\def\Arc{\frak{Arc}}

\def\sd{{\rm Sd}}
\def\sdd{{\rm Sd}^*}
\def\dop{{\Delta^{\rm op}}}

\def\lbt{{\tilde \Lambda}}

\def\ps{\id}
\def\pee{{\rm Perm}}
\def\epi{{(\tilde\Lambda^{\rm op})^\wedge}}
\def\lbto{{\lbt^{\rm op}}}
\def\nt{\N^{\times}}
\def\wnt{{\widehat{\N^{\times}}}}
\def\frg{\frak{g}}
\def\zmax{{\Z_{\rm max}}}

\def\totgs{oriented groupo\"ids }
\def\totg{oriented groupo\"id }
\def\arcnt{\Arc\ltimes\N^\times}
\begin{document}

\title{The Cyclic and Epicyclic  Sites}
\author{Alain Connes}
\email{alain@connes.org}
\address{Coll\`ege de France,
3 rue d'Ulm, Paris F-75005 France\newline
I.H.E.S. and Ohio State University.}
\author{Caterina Consani}
\email{kc@math.jhu.edu}
\address{Department of Mathematics, The Johns Hopkins
University\newline Baltimore, MD 21218 USA.}
%
%\dedication{}
\thanks{The second author 
%is partially supported by the NSF grant DMS 1069218 and 
would like to thank the Coll\`ege de France for some financial support.}
%This is the abstract\newline (see \textbf{elsdoc.pdf} on the svn and/or one of the sample %files
%\code{cmguide1.pdf} or \code{cmguide1.ps} available at
%\url{http://www.compositio.nl/cmauthor.html}).

\begin{abstract}
We determine the points of the  epicyclic topos which plays a key role in the geometric encoding of cyclic 
homology and the lambda operations. We show that the category of points of the  epicyclic topos is equivalent  to projective geometry in characteristic one
over algebraic extensions of the infinite semifield of ``max-plus integers'' $\Z_{\rm max}$. 
An object of this category is a pair $(E,K)$ of a \mod~$E$ over an algebraic extension $K$ of $\Z_{\rm max}$. The morphisms are projective classes of semilinear maps between \mods. The epicyclic topos sits over the arithmetic topos $\widehat{\N^\times}$  of \cite{CCas} and the fibers of the associated geometric morphism correspond to the cyclic site. 
In two appendices we review the role of the cyclic and epicyclic toposes as  the geometric structures supporting cyclic
homology and the lambda operations.\end{abstract}

\maketitle
\vspace*{1pt}
\tableofcontents  % for this guide only.
% A table of contents should normally not be included

%\newpage

\section{Introduction}
The theory of topoi of Grothendieck provides the best geometric framework to understand cyclic homology and the $\lambda$-operations using the topos associated to the cyclic category \cite{CoExt} and its epicyclic refinement \cite{CCproj}. Given a small category $\cC$, we denote by $\hat \cC$ the topos of contravariant functors from $\cC$ to the category of sets $\Se$. The epicyclic topos $\epi$ is obtained by taking the {\em opposite} of the epicyclic category $\lbt$. This choice   is dictated by the following natural construction. A commutative ring $R$ determines a  {\em covariant} functor $\fin\longrightarrow\Ab$ from the category of finite sets to that of abelian groups. This functor assigns to a finite set $J$ the tensor power $R^{\otimes J}=\bigotimes_{j\in J}R$. As explained in geometric terms here below, there is also a natural {\em covariant} functor  $\lbt\longrightarrow\fin$. The composite of these two functors $\lbt\longrightarrow\Ab$ provides, for any  commutative ring $R$,   a {\em covariant} functor $R^\natural$ from the epicyclic category to the category of abelian groups. In geometric terms $R^\natural$ is a sheaf of abelian groups over the topos $\epi$.  Both the cyclic homology of $R$ and its $\lambda$-operations are completely encoded by the associated sheaf $R^\natural$. In \cite{CCproj}, we provided a conceptual understanding of the epicyclic category as projective geometry over the semifield $\F:=\zmax$ of the tropical integers. In these terms the functor $\lbt\longrightarrow \fin$ considered above assigns to a projective space the underlying finite set. This article pursues the relation between the epicyclic topos and (projective) geometry in characteristic $1$ in more details. Our main result is the following (\cf~Theorem \ref{main})
\begin{theorem}\label{main} The category of points of the epicyclic topos $\epi$ is equivalent to the 
category $\cP$ whose objects are  pairs $(K,E)$, where $K$ is an algebraic extension of $\F=\Z_{\rm max}$ and $E$ is an archimedean \mod~over $K$. The  morphisms in $\cP$ are projective classes of semilinear maps and injective semifield morphisms. 
\end{theorem}

The ambigu\"ity in the choice of a representative of a projective class of semilinear maps in $\cP$  is inconvenient  when working, for example, with colimits. In Section \ref{sectN} (\S\ref{secttogs}) we provide a description of  the cyclic and the epicyclic categories in terms of a category $\mathfrak g$  of {\em \totgs} whose morphisms are no longer given by equivalence classes. There are by now a number of equivalent descriptions of the cyclic and epicyclic categories, ranging from the most concrete \ie given in terms of generators and relations, to the most conceptual as in \cite{CCproj}. The description of these categories  in terms of \totgs turns out to be very useful to determine the points of the  epicyclic topos by considering filtering colimits, in the category $\mathfrak g$, of the special points provided by  the Yoneda embedding of the categories. It is in fact well known that any point of a topos of the form $\hat \cC$  is obtained as a filtering colimit, in the category of flat functors $\cC\longrightarrow \Se$, of these special points. On the other hand, there is no guarantee ``a priori" that   this colimit process   yields the same result as the colimit  taken in the category $\mathfrak g$. This matter is solved  in two steps and in concrete terms  in Section \ref{sectF1}. In Proposition \ref{propmain} we  show how to associate to a pair $(K,E)$ as in the above Theorem a point of $\epi$. Conversely, in  \S\S\ref{sectmod}-\ref{recsect}  we explain a geometric procedure that allows one to  reconstruct the structure of an \totg from the flat functor naturally associated to a point of $\epi$.\newline
In \S \ref{sectgeomrel} we explore the relations  of $\epi$ with the arithmetic site $\wnt$, as recently defined in \cite{CCas}. Let $\nt$ be the small category  with a single object $\bullet$ and whose endomorphisms $\End(\bullet)=\nt$ form the multiplicative semigroup $\nt$ of  positive  integers. One has a canonical functor $\Mod:\lbto\longrightarrow \nt$ which is trivial on the objects and associates to a semilinear map of semimodules over $\F=\Z_{\rm max}$ the corresponding injective endomorphism $\fr_n\in \End(\F)$ (\cf \cite{CCproj} for details).   This functor  induces a geometric morphism of topoi $\Mod: \epi\longrightarrow \wnt$. The subcategory of $\lbto$ which is the kernel of this morphism  is the cyclic category $\Lambda$ ($\Lambda \simeq \Lambda^{\rm{op}}$).\newline
In the two Appendices we  interpret some known results on the simplicial, cyclic and epicyclic categories  in geometric terms by using the language of topos theory. It is important to realize the relevance of the language of Grothendieck topoi to interpret, for instance, the action of the barycentric subdivision on the points of the simplicial topos $\hat\Delta$. This language is also crucial  in order to carefully identify the natural generalization of cyclic homology and $\lambda$-operations in the framework of epicyclic modules. Appendix 
 \ref{appcyclic} stresses the nuance between $\lbt$ and $\lbto$ in a hopefully clear form.\newline
Appendix \ref{appbary} is dedicated to the description of the action of the barycentric subdivision on the points of the topos $\hat\Delta$. It is well known (\cf \cite{MM}) that these points correspond to intervals, \ie totally ordered sets $I$ with a smallest element $b$ and a largest element $t\neq b$. For each integer $k>0$, the barycentric subdivision $\sd_k$ defines an endofunctor of the simplicial category
$\Delta$  and one obtains in this way an action of the mono\"id $\nt$ by geometric morphisms on the topos $\hat\Delta$. We show that  the action of the barycentric subdivision on the points of $\hat\Delta$ is given by the operation of concatenation of $k$ copies of the interval $I$: the intermediate top point $t_j$ of the copy $I_j$ is identified with the bottom point $b_{j+1}$ of the subsequent copy $I_{j+1}$. Then, we form the small category $\dop\ltimes \N^\times$ crossed product of $\dop$  by the transposed action $\sdd$ of $\nt$ (\ie $\sd_{k}(f)^*=\sdd_k(f^*)$, where $f\mapsto f^*$ is the anti-isomorphism $\Delta\longrightarrow\dop$). This process allows one to view the $\lambda$-operations as 
elements $\Lambda_n^k$  of the associated convolution ring $\Z[\dop\ltimes \N^\times]$ with integral coefficients. We finally review  the geometric meaning of the $\lambda$-operations and the geometric proof of their commutation (\cf \cite{MCarthy}) with the Hochschild boundary operator. \newline
Appendix \ref{appcyclic} is dedicated to the description of the cyclic homology of cyclic modules (\cf \cite{CoExt}) and its extension to epicyclic modules \cite{Loday}. An  {\em epicyclic} module $E$ is a  {\em covariant} functor  $\lbt\longrightarrow\Ab$. These modules correspond to sheaves of abelian groups on the topos $\epi$. At this point the nuance between the epicyclic category and its dual  plays an important role since unlike the cyclic category the epicyclic category is not anti-isomorphic to itself. As explained earlier on in this introduction, a commutative ring $R$ gives rise naturally to  an epicyclic module $R^\natural$ and it  is well known (\cf \cite{Loday}) that the $\lambda$-operations on cyclic homology of $R$ are obtained directly through the associated epicyclic module. We provide a simple and conceptual proof of the commutation of the $\lambda$-operations with the $B$ operator of cyclic theory. Finally, we point out that the extended framework of epicyclic modules involves many more modules than those arising by composition, as explained earlier, from a covariant functor $\fin\longrightarrow\Ab$. In fact, these  particular (epicyclic) modules have {\em integral} weights and the $\lambda$-operations decompose their cyclic homology as direct sums of modules on which $\Lambda_n^k$ acts by an integral power of $k$. This integrality property no longer holds for general epicyclic modules as  can be easily checked by applying a twisting argument.

%Section: The epicyclic category and oriented groupoids
 
\section{The epicyclic category and the \totgs}\label{sectN}

%Subsection: Generalities on oriented groupoids
\subsection{Generalities on \totgs}
A groupo\"id $G$ is a small category where  the morphisms are invertible. Given a subset $X\subset G$ of a groupo\"id, we set $X^{-1}:=\{\gamma^{-1}\mid \gamma\in X\}$. Let $G^{(0)}$ be the set of objects of $G$ and denote by
$r,s:G\to G^{(0)}$ the range and the source maps respectively. We view $G^{(0)}$ as the subset  of units of $G$. The following  definition is a direct generalization to groupo\"ids of the notion of right ordered group (\cf\cite{glass})
\begin{defn}  An \totg  $(G,G_+)$ is a groupo\"id $G$ endowed with a subcategory $G_+\subset G$, such that the following relations hold
\begin{equation}\label{orddefn}
G_+\cap G_+^{-1}=G^{(0)}, \qquad G_+\cup G_+^{-1}=G.
\end{equation}
\end{defn}

Let  $(G,G_+)$ be an \totg and  
let $x\in G^{(0)}$. The set $G_x:=\{\gamma\in G\mid s(\gamma)=x\}$ is endowed with the total order defined by
\begin{equation}\label{iso8}
\gamma\leq \gamma'\iff \gamma'\circ \gamma^{-1}\in G_+.
\end{equation}
This order is right invariant by construction: \ie  for any $\beta\in G$, with $r(\beta)=x$, one has
$$
\gamma\leq \gamma'\iff \gamma\circ \beta\leq \gamma'\circ \beta.
$$
In the following subsections we describe two constructions of \totgs associated to a group action. 

%Subsubsection
\subsubsection{$G=X\ltimes H$.}\label{subsubsemidirect}

Let $H$ be a group acting on a set $X$. Then the semi-direct product $G:=X\ltimes H$ is a groupo\"id with source, range and composition law defined respectively as follows
$$
s(x,h):= x, \quad r(x,h):= hx, \quad (x,h)\circ (y,k):=(y,hk). 
$$
(As in any groupo\"id the composition $\gamma\circ \gamma'$ is only defined when $s(\gamma)=r(\gamma')$ which holds here if and only if  $x=ky$).
One has a canonical homomorphism of groupo\"ids $\rho:G\to H$, $\rho(x,h)=h$.
\begin{lem} \label{examsemi0} Let $(H,H_+)$ be a right ordered group. Assume that  $H$ acts on a set $X$. Then the semi-direct product $G=X\ltimes H$ with $G_+:=\rho^{-1}(H_+)$ is an \totg\!\!.
\end{lem}
\proof
By definition, the subset $H_+\subset H$ of the group $H$ is stable under product and fulfills the equalities: $H_+\cap H_+^{-1}=\{1\}, \    \  H_+\cup H_+^{-1}=H$. This implies \eqref{orddefn} using $\rho^{-1}(\{1\})=G^{(0)}$. \qed

Let, in particular, $(H,H_+)=(\Z,\Z_+)$ act by translation on the set
$X=\Z/(m+1)\Z$ of integers modulo $m+1$. Then one obtains the \totg 
\begin{equation}\label{examsemi}
\frg(m):=(\Z/(m+1)\Z)\ltimes \Z.
\end{equation}
The
\totgs $\frg(m)$ will play a crucial role in this article.

%\subsubsection
\subsubsection{$G=(X\times X)/H$.}\label{subsubsquare}

Let $H$ be a group acting freely on a set $X$. 
Let   $G(X,H)=(X\times X)/H$ be the quotient of $X\times X$ by the diagonal action of $H$
$$
G(X,H):=(X\times X)/\sim \qquad  (x,y)\sim (h(x),h(y)), \  \forall h\in H.
$$
Let $r$ and  $s$ be the two projections of $G(X,H)$ on $G^{(0)}:=X/H$ defined by $r(x,y)=x$ and $s(x,y)=y$. Let $\gamma, \gamma'\in G(X,H)$ be such that $s(\gamma)=r(\gamma')$. Then, for $\gamma\sim (x,y)$ and $\gamma'\sim (x',y')$, there exists a unique  $h\in H$ satisfying $x'=h(y)$: this because  $s(\gamma)=r(\gamma')$ and $H$ acts freely on $X$. Then, the pair $(h(x),y')$ defines an element of $G(X,H)$ independent of the choice of the pairs  representing the elements $\gamma$ and $\gamma'$. We denote by $\gamma\circ \gamma'$ the class of $(h(x),y')$ in $G(X,H)$. This construction defines a groupo\" id law on $G(X,H)=(X\times X)/H$.

\begin{lem} \label{examsq} Let $H$ be a group acting freely on a set $X$.  Assume that  $X$ is totally ordered and that $H$ acts by order automorphisms. Then $G(X,H)=(X\times X)/H$ is an \totg\!\! with 
\begin{equation}\label{ordgrp}
G_+(X,H)=\{(x,y)\in G(X,H)\mid x\geq y\}.
\end{equation}
\end{lem}
\proof Since $H$ acts by order automorphisms the condition $x\geq y$ is independent of the choice of a representative $(x,y)$ of a given $\gamma\in G(X,H)=(X\times X)/H$. This condition defines a subcategory $G_+(X,H)$ of $G(X,H)$. The conditions \eqref{orddefn} then follow since $X$ is totally ordered.\qed

\begin{lem}\label{subsubfrg}
Let $X=\Z$ with the usual total order. Let $m\in \N$ and let the group $\Z$ act on $X$ by $h(x):=x+(m+1)h$, $\forall x\in X, h\in\Z$. Then the \totg $G=(X\times X)/\Z$ is canonically isomorphic to the \totg $\frg(m)$ of \eqref{examsemi}.
\end{lem}
\proof  The associated \totg $(G,G_+)$ is by construction the quotient of $\Z\times \Z$ by the equivalence relation: $(x,y)\sim (x+\ell(m+1),y+\ell(m+1))$, $\forall \ell\in \Z$. Thus the following map defines a bijective homomorphism of groupoids 
$$
\psi: G\to \frg(m)=(\Z/(m+1)\Z)\ltimes \Z, \qquad \psi(x,y)=(\pi(y),x-y)
$$
where $\pi:\Z\to \Z/(m+1)\Z$ is the natural projection. 
One has by restriction $G_+\stackbin[\sim]{\psi}{\to} \frg_+(m)$, since $x\geq y\iff x-y\geq 0$, so that $\psi$ is in fact an isomorphism of \totgs\!\!.\qed

%Subsection: The oriented groupoid associated to an archimedean set

\subsection{The \totg associated to an archimedean set}\label{sectarchtogr}

In this section we explain how to associate an \totg to an archimedean set and describe the special properties of the \totgs thus obtained. We first  recall from \cite{topos} the definition of an archimedean set.
\begin{defn}\label{defnarc} An archimedean set is a pair $(X,\theta)$ of
a non-empty, totally ordered set $X$  and an order automorphism $\theta\in \Aut X$, such that $\theta(x)>x$, $\forall x\in X$. The automorphism $\theta$ is also required to fulfill the following archimedean property
\begin{equation*}%\label{arch}
\forall x,y\in X, \ \exists n\in \N~{\rm s.t.} \quad y\leq \theta^n(x).
\end{equation*}
\end{defn}

Let  $(X,\theta)$ be an archimedean set and let $G(X,\theta)$ be the \totg associated by Lemma \ref{examsq} to the action of $\Z$ on $X$ by integral    powers of $\theta$. Thus
$$
G(X,\theta):=(X\times X)/\sim, \  \  (x,y)\sim (\theta^n(x),\theta^n(y)), \  \forall n\in \Z
$$
and
\begin{equation}\label{ordgrpbis}
G_+(X,\theta):=\{(x,y)\in G(X,\theta)\mid x\geq y\}.
\end{equation}
Next proposition describes the properties of the pair $(G(X,\theta),G_+(X,\theta))$ so obtained.
\begin{prop} \label{proparctog} The \totg $(G,G_+)=(G(X,\theta),G_+(X,\theta))$  fulfills the following conditions
\begin{enumerate}
\item  $\forall x,y\in G^{(0)}$, $\exists \gamma\in G_+$~s.t.~ $s(\gamma)=y$, 
$r(\gamma)=x$.

\item For $x\in X$, the ordered groups $G_x^x:=\{\gamma\mid s(\gamma)=r(\gamma)=x\}$ are isomorphic to  $(\Z,\leq)$.

\item Let $\gamma\in G$ with $s(\gamma)=y$ and  
$r(\gamma)=x$. Then the map: $G_y^y\ni\rho\mapsto \gamma\circ \rho\circ \gamma^{-1}\in G_x^x$ is an isomorphism of ordered groups.
\end{enumerate}

\end{prop} 
\proof Since $\theta$ is an order automorphism of $X$, the group $\Z$ acts by order automorphisms. We check the three conditions (i)-(iii). 

$(i)$~For $x,y\in X$, there exists $n\in \N$ such that $\theta^n(x)\geq y$. Then $\gamma\sim (\theta^n(x),y)$ belongs to $G_+$ and $s(\gamma)= y$, $r(\gamma)= x$.

$(ii)$~Let $x\in X$. The conditions $s(\gamma)=r(\gamma)= x$ imply that the class of $\gamma\in G(X,\theta)$ admits a unique representative of the form $(\theta^n(x),x)$. One easily checks that the map $G_x^x \to (\Z,\leq)$,  $(\theta^n(x),x)\mapsto n$ is an isomorphism of ordered groups.

$(iii)$~Let $x,y\in X$ with $\gamma\sim (x,y)$. Then for $\rho\sim (\theta^n(y),y)\in G_y^y$ one gets $ \gamma\circ \rho\circ \gamma^{-1}\sim (\theta^n(x),x)$,  thus the unique isomorphism with $(\Z,\leq)$ is preserved. \qed

Let  $(G,G_+)$ be an \totg fulfilling the three conditions of Proposition \ref{proparctog}. 
Let 
$x\in G^{(0)}$, consider the set $G_x=\{\gamma\in G\mid s(\gamma)=x\}$  with the total order defined by \eqref{iso8}
and with the action of $\Z$ given, for $\gamma_x \in G_+$ the positive generator of $G_x^x$, by
\begin{equation}\label{iso7}
\theta(\gamma):=\gamma\circ \gamma_x. 
\end{equation}
When one applies this construction to the case $(G,G_+)=(G(X,\theta),G(X,\theta)_+)$, for $(X,\theta)$ an archimedean set, with $x\in G^{(0)}=X/\theta$ one obtains, after choosing a lift $\tilde x\in X$ of $x$, an isomorphism 
\begin{equation}\label{iso9}
j_{\tilde x}:X \stackrel{\sim}{\to} G_x, \qquad  j_{\tilde x}(z)=(z,x).
\end{equation}

The following proposition shows that the two  constructions $(X,\theta)\mapsto G(X,\theta)$ and $(G,G_+)\mapsto (G_x,\theta)$  are  reciprocal.  
\begin{prop} \label{proparctogbis} Let  $(G,G_+)$ be an \totg fulfilling the conditions of Proposition \ref{proparctog} and  
let 
$x\in G^{(0)}$. Consider the set $G_x:=\{\gamma\in G\mid s(\gamma)=x\}=X$ endowed with the total order \eqref{iso8} and the action of $\Z$ on it given by \eqref{iso7}. Then $(X,\theta)$ is an archimedean set and one has an isomorphism of  \totgs  $$(G(X,\theta),G(X,\theta)_+)\cong (G,G_+).$$
\end{prop}

\proof The implication 
 $\gamma\leq \gamma'\implies \theta(\gamma)\leq \theta( \gamma')$ follows since right multiplication preserves the order. Moreover, the condition $(iii)$ of Proposition \ref{proparctog} implies that  $\theta(\gamma)>\gamma$, $\forall \gamma\in X=G_x$. Next we show that the archimedean property holds on $(X,\theta)$. Let  $\gamma\leq \gamma' \in G_x$, with 
$y=r(\gamma)$ and $y'=r(\gamma')$. By applying the condition $(i)$ of Proposition \ref{proparctog}, we choose $\delta\in G_+$ such that $s(\delta)=y'$ and $r(\delta)=y$. Then $\gamma''=\delta\circ \gamma'$  fulfills $s(\gamma'')=s(\gamma)$ and $r(\gamma'')=r(\gamma)$ and thus there exists $n\in \Z$ such that $\gamma''=\gamma\circ \gamma_x^n$. Moreover, one has $\gamma''=\delta\circ \gamma'\geq \gamma'$.
It follows that  $(X,\theta)$ is an archimedean set. If one replaces $x\in G^{(0)}$ by $y\in G^{(0)}$, then the condition $(i)$ implies that there exists  $\alpha\in G_+$ with $s(\alpha)=y$, 
$r(\alpha)=x$. Then the map $G_x\to G_y$, $\gamma\mapsto \gamma\circ \alpha$ is an order isomorphism which satisfies 
$$
(\gamma\circ \gamma_x)\circ \alpha=(\gamma\circ \alpha)\circ (\alpha^{-1}\circ  \gamma_x\circ \alpha).
$$
Since condition $(iii)$ implies $\alpha^{-1}\circ  \gamma_x\circ \alpha=\gamma_y$, one obtains an isomorphism of the corresponding archimedean sets.

Finally, we compare the pair $(G,G_+)$ with $(G(X,\theta), G_+(X,\theta))$. We define a map $f:G(X,\theta)\to G$ as follows: given a pair  $(\gamma,\gamma')$ of elements of $X=G_x$, one sets $f(\gamma,\gamma'):=\gamma\circ \gamma'^{-1}$. One has 
$$
f(\theta(\gamma),\theta(\gamma'))=f(\gamma\circ \gamma_x,\gamma'\circ \gamma_x)=\gamma\circ \gamma'^{-1}=f(\gamma,\gamma').
$$
To show that $f$ is a groupo\"id homomorphism it is enough to check that $f(\gamma,\gamma')\circ f(\gamma',\gamma'')=f(\gamma,\gamma'')$ and this can be easily verified. Next we prove that $f$ is bijective. Let $\alpha\in G$. By  applying condition $(i)$ of Proposition \ref{proparctog}, there exists $\gamma\in G$ such that $r(\gamma)=r(\alpha)$ and $s(\gamma)=x$. Let then $\gamma'=\alpha^{-1}\gamma$. Since $s(\gamma')=x$, both $\gamma,\gamma'$ belong to $X=G_x$ and moreover  $f(\gamma,\gamma')=\alpha$  showing that $f$ is surjective.  Let $\gamma_j,\gamma'_j$ be elements of $X=G_x$ such that 
$f(\gamma_1,\gamma'_1)=f(\gamma_2,\gamma'_2)$. One then has $\gamma_1\circ \gamma'^{-1}_1=\gamma_2\circ \gamma'^{-1}_2$ and hence $\gamma^{-1}_2\circ \gamma_1=\gamma'^{-1}_2\circ \gamma_2=\gamma_x^n$ for some $n\in \Z$. It follows that $(\gamma_	1,\gamma'_1)=(\theta^n(\gamma_2),\theta^n(\gamma'_2))$ which shows that $f$ is also injective. Finally, 
for any $\gamma,\gamma'\in X=G_x$ one has
$
\gamma\geq \gamma'\iff \gamma\circ \gamma'^{-1}\in G_+
$
showing that $f$, so defined, is an order isomorphism.\qed

%The category of archimedean sets in terms of oriented groupoids

\subsection{The category of archimedean sets in terms of \totgs}\label{secttogs}
In this section we extend the above construction $(X,\theta)\to G(X,\theta)$ of the oriented groupoid associated to an archimedean set to a functor $G$ connecting the category $\arcnt$ of archimedean sets to that of \totgs\!\!. 
We recall, from \cite{topos}, the definition of the category of archimedean sets.

\begin{defn}\label{last2} The objects of the category $\arcnt$ are the archimedean sets $(X,\theta)$ as in Definition \ref{defnarc},  the morphisms  
$f: (X,\theta)\to (X',\theta')$ in $\arcnt$ are  equivalence classes of maps 
\begin{equation}\label{deff2}
f:X\to X', \ \ f(x)\geq f(y) \  \  \forall x\geq y ;\qquad \exists k>0, \  f(\theta(x))=\theta'^k(f(x)),\quad \forall x \in X
\end{equation}
where the equivalence relation identifies two such maps $f$ and $g$ if there exists an integer $m\in \Z$ such that 
$g(x)=\theta'^m(f(x))$, $\forall x\in X$.
\end{defn}
Given a morphism  
$f: (X,\theta)\to (X',\theta')$ in $\Arc\ltimes \N^\times$ as in \eqref{deff2}, the map $G(f)$ sending $(x,y)\mapsto (f(x),f(y))$ is well defined  since
$$
X\times X \ni (x,y)\sim (x',y')~\Longrightarrow~ (f(x),f(y))\sim (f(x' ),f(y')).
$$
Moreover  $G(f):G(X,\theta)\to G(X',\theta')$ does not change if we replace $f$ by $g(x)=\theta'^m(f(x))$ since this variation does not alter the element $(f(x),f(y))\sim (g(x),g(y))\in G(X',\theta')$.

\begin{prop}\label{propfaith} The association $(X,\theta)\mapsto G(X,\theta)$, $f\mapsto G(f)$ defines a faithful functor $G$ from $\Arc\ltimes \N^\times$ to the category of oriented groupoids. For any non-trivial morphism of \totgs $\rho: G(X,\theta)\to G(X',\theta')$ there exists  a (unique) morphism  $$f\in \Hom_{\arcnt}((X,\theta),(X',\theta')) \quad{\rm s.t.}\quad  \rho=G(f).$$
\end{prop}
\proof  By construction the map $G(f)$ is a morphism of groupo\"ids since 
$$
G(f)\left((x,y)\circ (y,z)\right)=(f(x),f(z))=G(f)((x,y))\circ G(f)((y,z))
$$
Moreover it is a morphism of \totgs since by \eqref{deff2}, one has $f(x)\geq f(y)$, $ \forall x\geq y$ so that 
$$
(x,y)\in G_+(X,\theta)\implies (f(x),f(y))\in G_+(X,\theta)
$$
We show that the functor $G$ is faithful. Assume that  $f,g\in \Hom_{\arcnt}((X,\theta),(X',\theta'))$ are such that $G(f)=G(g)$. Then for any $(x,y)\in X\times X$ there exists an integer $n=n(x,y)\in \Z$ such that 
$$
g(x)=\theta'^{n(x,y)}(f(x)), \  \  g(y)=\theta'^{n(x,y)}(f(y)).
$$
Since $\theta'$ acts freely on $X'$, the integer $n(x,y)$ is unique. The first and the second equations prove that $n(x,y)$ is independent of $y$ and $x$ respectively. Thus one derives that $f$ and $g$ are in the same equivalence class, \ie they define the same element of $\Hom_{\arcnt}((X,\theta),(X',\theta'))$.

Let now $\rho: G(X,\theta)\to G(X',\theta')$ be a non-trivial morphism of \totgs\!\!. Let $x\in G^{(0)}$, $y=\rho(x)$: the ordered group morphism  $\rho:G_x^x\to G'^y_y$ is non-constant and is given by $\rho(\gamma_x)=\gamma_y^k$, for some $k>0$. The map $\rho:G_x\to G'_y$ is non-decreasing since
$$
\gamma\leq \gamma'\implies \gamma'\circ \gamma^{-1}\in G_+\implies \rho(\gamma')\circ \rho(\gamma)^{-1}\in G'_+.
$$
Given an archimedean set $(X,\theta)$ and an element $z\in X$,   we claim that the map $\psi_{X,z}:X\to G(X,\theta)$ defined by 
$
\psi_{X,z}(y):=(y,z)
$
is an order preserving bijection of $X$ with $G(X,\theta)_x$,  where $x$ is the class of $z$ in $X/\theta$. Indeed, every element of $G(X,\theta)_x$ admits a unique representative of the form $(y,z)$ and one has
$$
\psi_{X,z}(y)\leq \psi_{X,z}(y')\iff (y',z)\circ (z,y) \in G_+\iff y'\geq y.
$$
Moreover $\psi_{X,z}$ is also equivariant since one has
$$
\psi_{X,z}(\theta(y))=(\theta(y),z)=(\theta(y),\theta(z))\circ (\theta(z),z)=\psi_{X,z}(y)\circ \gamma_x.
$$
Let then $\tilde x\in X$ and $\tilde y \in X'$ be two lifts of $x$ and $y$. The  map  
$
f:=\psi_{X',\tilde y}^{-1}\circ \rho\circ \psi_{X,\tilde x}
$
is a non-decreasing map from $X$ to $X'$ and one has, using $\rho(\gamma\circ \gamma_x)=\rho(\gamma)\circ \gamma_y^k$, that
\begin{equation*}
  f(\theta(x))=\theta'^k(f(x)),\  \forall x \in X.
\end{equation*}
One derives by construction the equality $(f(a),\tilde y)=\rho((a,\tilde x)$, $\forall a\in X$. Taking $a=\tilde x$ this gives $f(\tilde x)=\tilde y$ since $\rho((\tilde x,\tilde x))$ is a unit. This shows that $\rho(\gamma)=G(f)(\gamma)$ $\forall\gamma\in G_x$ and the same equality holds for all $\gamma \in G$ since both $\rho$ and $G(f)$ are homomorphisms while any element of $G$ is of the form $\gamma\circ (\gamma')^{-1}$ with $\gamma, \gamma'\in G_x$.\qed

\begin{rem}\label{epicor0} {\rm
The \totgs  associated to archimedean sets are all equivalent, in the sense of equivalence of (small) categories, to the ordered group $(\Z,\Z_+)$. It follows that a morphism $\phi$  of \totgs induces an associated morphism $\Mod(\phi)$ of totally ordered groups, \ie an ordered group morphism $(\Z,\Z_+)\to (\Z,\Z_+)$ given by multiplication by an integer $\Mod(\phi)=k\in \N$. 
Proposition \ref{propfaith} suggests to refine the category $\frg$ of \totgs  by considering only the morphisms $\phi$ such that $\Mod(\phi)\neq 0$. In other words what one requires is that the associated morphism of totally ordered groups, obtained by working modulo equivalence of categories, is {\em injective}. One can then reformulate Proposition \ref{propfaith}  stating that the functor $G$ is full and faithful.}
\end{rem}

\begin{cor}\label{epicor} The epicyclic category $\lbt$ (\cf~Appendix~\ref{appcyclic}) is canonically isomorphic to the category with objects the \totgs $\frg(m)$, $m\geq 0$, of Example \ref{examsemi} and morphisms the non-trivial morphisms of \totgs\!\!. \newline The functor which  associates to a morphism of \totgs its class up to equivalence coincides with the functor $\Mod:\lbt\longrightarrow \nt$ which sends a semilinear map of semimodules over $\F=\Z_{\rm max}$ to the corresponding injective endomorphism $\fr_n\in \End(\F)$ (\cf \cite{CCproj})
\end{cor}
\proof By Proposition 2.8 of \cite{CCproj}, the epicyclic category $\lbt$ is canonically isomorphic to the full subcategory of $\Arc\ltimes\N^\times$ whose objects are the archimedean sets $\underline m:=(\Z,x\mapsto x+m+1)$ for $m\geq 0$. The \totg $G(\underline m)$ is by 
Lemma \ref{subsubfrg} canonically isomorphic to $\frg(m)$. The first statement then follows from Proposition \ref{propfaith} while the last one is checked easily and directly.\qed

%Section: Points of the epicyclic topos and projective geometry in characteristic one

\section[Points of the epicyclic topos $\epi$]{Points of $\epi$ and  projective geometry in characteristic one}\label{sectF1}

Let $\cC$ be a small category and $\hat \cC$ the topos of contravariant functors from $\cC$ to the category of sets $\Se$. Yoneda's Lemma defines  an embedding of the opposite category $\cC^{\rm op}$ into the category of points  of the topos $\hat \cC$. More precisely, to an object $c$ of $\cC^{\rm op}$ one associates the flat covariant functor $h_c(-):\cC\longrightarrow\Se$, $h_c(-)=\Hom_\cC(c,-)$. Then, one sees that through  Yoneda's embedding $h: \cC^{\rm op}\longrightarrow \hat\cC$, $c\mapsto h_c$,  any point of $\hat \cC$ can be obtained as a filtering colimit of points of the form $h_c$. We shall apply these well known general facts to $\cC=\lbto$: the opposite of the epicyclic category $\lbt$. We refer to Appendix~\ref{appcyclic} for the basic notations on the cyclic and epicyclic categories. It follows from  Corollary \ref{epicor}  that  $\lbt$ is canonically isomorphic to the full subcategory   of the category $\frg$ of \totgs whose objects are the \totgs of the form $\frg(m)$. This fact suggests that one should obtain the points of the topos $\epi$ by considering filtering colimits of the objects $\frg(m)$ in  $\frg$. In this section we compare the colimit procedures taken respectively in the category of flat functors $\lbto\longrightarrow\Se$ and    in the category $\frg$. The comparison is made directly by reconstructing the structure of an \totg starting from a flat functor as above.
The main result  is the following
\begin{thm}\label{main} The category of points of the epicyclic topos $\epi$ is equivalent to the 
category $\cP$ whose objects are  pairs $(K,E)$ where $K$ is an algebraic extension of $\F=\Z_{\rm max}$ and $E$ is an archimedean \mod~over $K$. The  morphisms are projective classes of semilinear maps and injective semifield morphisms. 
\end{thm}
One knows from \cite{CCas} that an algebraic extension $K$ of the semifield $\F=\Z_{\rm max}$ of tropical integers is  equivalently described by a totally ordered group $(H,H_+)$ isomorphic to a subgroup $\Z\subseteq H\subseteq \Q$ of the rationals. 
An archimedean \mod~$E$ over $K$ is in turn  described  (\cf\cite{CCproj}) by a totally ordered set $X$ on which $H$ acts by order automorphisms of type: $(x,h)\mapsto x+h$ which fulfill the property
\begin{equation}\label{free}
h+x>x, \quad \forall h\in H_+,~ h\neq 0, \  x\in X
\end{equation}
and the archimedean condition
\begin{equation}\label{archH}
\forall x,y\in X, \ \exists h\in H_+ ~{\rm s.t.}\quad  h+x>y.
\end{equation}
It follows from \cite{AGV} (\cf~also \cite{MM}, Theorem 2 Chapter VII \S 5) that a point of a topos of the form $\hat\cC$, where $\cC$ is a small category, is described by a covariant {\em flat} functor $F: \cC\longrightarrow \Se$.
Next, we overview the strategy adopted to prove Theorem \ref{main}. 

In \S\ref{step1} we associate to a pair $(K,E)$ a point of $\epi$. This construction is accomplished in two steps. First, we extend the construction $(X,\theta)\mapsto G(X,\theta)$  of the \totg associated to an archimedean set (as in \S\ref{sectarchtogr})   to a pair $(K,E)$ as in Theorem \ref{main}. Then, for any given pair $(K,E)$, we provide a natural construction of a point of $\epi$ by means of the following associated {\em flat} functor ($\underline n =(\Z,x\mapsto x+n+1)$, $n\geq 0$)
\begin{equation}\label{assflat}
F:\lbto \longrightarrow \Se, \qquad   F(\underline n)=\Hom_\frg\left(\frg(n),G(K,E)\right).
\end{equation}
Here, one implements Corollary \ref{epicor} to identify the category $\lbt$ with a  full subcategory of the category $\frg$ of \totgs with injective morphisms (up to equivalence).

To produce the converse of the above construction, \ie in order to show that any point of $\epi$ is obtained as in  \eqref{assflat} by means of a uniquely associated pair $(K,E)$, we start from  a covariant flat functor $F: \lbto\longrightarrow \Se$
and describe in \S\S\ref{sectmod}-\ref{recsect}  a procedure that allows one to reconstruct the semifield $K$ by using  the natural geometric morphism of topoi associated to the functor $\Mod:\lbto\longrightarrow \nt$. The archimedean semimodule $E$ (totally ordered set) is then reconstructed by using a suitable restriction of $F$ to obtain intervals from points of the simplicial topos $\hat \Delta$.

%Subsection: The flat functor associated to a pair

\subsection{The flat functor $\lbto\longrightarrow\Se$ associated to a pair $(K,E)$}\label{step1}

Let $(H,H_+)$ be a totally ordered abelian group, denoted additively and  $X$ a totally ordered set on which the ordered group $H$ acts preserving the order and fulfilling  \eqref{free}. Let $(G(X,H),G_+(X,H))$ be the  \totg associated  to the pair $(X,H)$ by Lemma \ref{examsq}, thus one has
\begin{equation}\label{free1}
G(X,H):=(X\times X)/H, \  \   G_+(X,H):=\{(x,y)\mid x\geq y\}.
\end{equation}
 The next lemma is used to show that the  functor $F: \lbto\longrightarrow \Se$ naturally associated  to a pair $(K,E)$ is {\em filtering}.

\begin{figure}
\begin{center}
\includegraphics[scale=0.6]{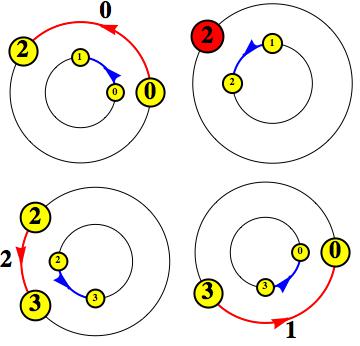}
\end{center}
\caption{One encodes a morphism $\phi$ of \totgs from $\frg(n)$ to $\frg(m)$ by the arrows $\gamma_i=\phi(\alpha(n,i))$ associated to the generators $\alpha(n,i)$. To each generator (in blue) one assigns an arrow (in red) specified by its source and range and by an integer which gives the number of additional windings.\label{clock1} }
\end{figure}

\begin{lem}\label{lemmain} 
Let $(H,H_+)$ be a non-trivial subgroup of $(\Q,\Q_+)$ and assume that the totally ordered set $X$ on which $H$ acts fulfills the  archimedean condition \eqref{archH}.
Let $\cF=\{\phi_j\mid 1\leq j\leq n\}$ be a finite set of morphisms $\phi_j\in \Hom_\frg(\frg(m_j),G(X,H))$.\newline
Then, there exists a singly generated subgroup $H_0\subset H$, a subset $X_0\subset X$ stable under the action of $H_0$, morphisms  $\psi_j\in \Hom_\frg(\frg(m_j),G(X_0,H_0))$  and an integer $m\in \N$ such that, denoting by $\iota:G(X_0,H_0)\to G(X,H)$ the natural morphism, one has 
\begin{equation}\label{archH0}
\phi_j=\iota\circ \psi_j\quad \forall j, \qquad  G(X_0,H_0)\simeq \frg(m).
\end{equation}
Moreover, let $\psi, \psi'\in \Hom_\frg(\frg(n),G(X_0,H_0))$ be two morphisms such that $\iota\circ \psi=\iota\circ \psi'$, then there exists a  singly generated subgroup $H_1$ with $H_0\subset H_1\subset H$, a subset $X_1\subset X$ containing $X_0$ and stable under the action of $H_1$, such that the equality $\iota_1\circ \psi=\iota_1\circ \psi'$ holds in $\Hom_\frg(\frg(n),G(X_1,H_1))$, with $\iota_1:G(X_0,H_0)\to G(X_1,H_1)$ the natural morphism.
\end{lem}
\proof  We denote by  
$
\alpha(m,i):=(i,1)\in \left(\Z/(m+1)\Z\right)\ltimes \Z=\frg(m)
$
the natural positive generators of the \totg $\frg(m)$. Let $(G,G_+)$ be an \totg\!\!. 
A morphism $\phi\in \Hom_\frg(\frg(m),G)$ is uniquely specified by the $m+1$ elements $\gamma_i=\phi(\alpha(m,i))\in G_+$ (\cf ~Figure~\ref{clock1}) fulfilling the conditions (with $\gamma_{m+1}:=\gamma_{0}$)
$$
r(\gamma_i)=s(\gamma_{i+1})\qqq i,~ 0\leq i \leq m, \quad \gamma_m\circ \cdots \circ \gamma_0\notin G^{(0)}.
$$
If $G=G(X,H)$, it follows that  $\Hom_\frg(\frg(m),G)$ is the quotient of the subset  of $X^{m+2}$
\begin{equation}\label{functn}\{(x_0,\ldots , x_{m+1})\in X^{m+2}\mid x_j\leq x_{j+1}, \forall j\leq n,~ x_{m+1}\in x_0+H_+\}
\end{equation}
by the diagonal action of $H$.
The morphism $\phi\in \Hom_\frg(\frg(m),G)$ associated to $(x_0,\ldots , x_{m+1})$ is given by 
\begin{equation}\label{functnbis}\phi(\alpha(m,i))= (x_{i+1},x_i)\in G_+\qqq i, ~ 0\leq i\leq m.
\end{equation}
For each $\phi_j\in \cF$ we get an $h_j\in H_+$, then we choose a finite subset $Z_j\subset X$ such that $\phi_j$ is represented by an (m+2)-uple $(x_0,\ldots , x_{m+1})\in X^{m+2}$ with all $x_i\in Z_j$. Let $H_0$ be the subgroup of $H$ generated by the $h_j$'s and let $X_0$ be the $H_0$ invariant subset of $X$ generated by the union of the $Z_j$. Then $H_0$ is singly generated and the pair $(X_0,H_0)$ fulfills the archimedean condition \eqref{archH}. Moreover by construction one can lift the maps $\phi_j\in \cF$ to elements $\psi_j\in \Hom_\frg(\frg(m_j),G(X_0,H_0))$ such that $\phi_j=\iota\circ \psi_j$, where $\iota:G(X_0,H_0)\to G(X,H)$ is the natural morphism.

It remains to show that $G(X_0,H_0)\simeq \frg(m)$ for some $m\in\N$. By construction there exists a finite subset $Z\subset X_0$ such that $X_0=Z+H_0$. Using the archimedean property \eqref{archH} and since $H_0\simeq \Z$, it follows that the pair $(X_0,H_0)$ is an archimedean set such that the quotient $X_0/H_0$ is finite. Thus  $G(X_0,H_0)\simeq \frg(m)$, where $m+1$ is the cardinality of $X_0/H_0$. 

To prove the last statement,  let $(x_0,\ldots , x_{m+1})\in X_0^{m+2}$ (resp. $(x'_0,\ldots , x'_{m+1})\in X_0^{m+2}$) represent $\psi$ (resp. $\psi'$). The equality $\iota\circ \psi=\iota\circ \psi'$ implies that there exists $h\in H$ such that $x'_j=x_j+h$ for all $j$. One then lets $H_1$ be the subgroup of $H$ generated by $H_0$ and $h$ and $X_1$ the $H_1$ invariant subset of $X$ generated by  $X_0$.\qed

\begin{prop}\label{propmain} 
Let $(H,H_+)$ be a non-trivial subgroup of $(\Q,\Q_+)$ and assume that the totally ordered set $X$ on which $H$ acts fulfills the  archimedean condition
\eqref{archH}.
Then the following formula defines a flat  functor 
\begin{equation}\label{archH1}
F: \lbto\longrightarrow\Se\qquad F(\underline n)=\Hom_\frg(\frg(n),G), \quad G=G(X,H)
\end{equation}
where $\underline n=(\Z, x\stackrel{\theta}{\mapsto}x+n+1)$.
\end{prop}
\proof The statement follows from Lemma \ref{lemmain} showing that the functor $F$  is obtained as a filtering colimit of {\em representable, flat} functors. We provide the detailed proof for completeness and to review the basic properties of flat functors which will be used later in this article. Corollary \ref{epicor} provides a canonical identification of the epicyclic category $\lbt$ with the full subcategory of $\frg$ of oriented groupoids of the form $\frg(m)$. In particular, \eqref{archH1} defines a covariant functor. It remains  to show that this functor is flat.
One knows from classical facts in the theory of Grothendieck topoi (\cf~\eg \cite{MM}, Chapter VII \S 6, Theorem 3) that a functor $F: \cC\longrightarrow \Se$ ($\cC$ a small category) is flat if and only if it  is filtering \ie  the category $\int_{\cC}\, F$ is filtering (\cf \cite{MM} Chapter VII \S 6, Definition 2).  The objects of the category $\int_{\cC}\, F$ are pairs $(j,x)$ where $j$ is an object of $\cC$ and $x\in F(j)$. The morphisms between two such objects $(j,x)$ and $(k,y)$ are elements  $\gamma\in\Hom_\cC(j,k)$ such that $F(\gamma)x=y$. We recall that the filtering condition on a small category $I$  is equivalent to the fulfillment of the following conditions

\begin{enumerate}
\item $I$ is non empty.
\item For any two objects $i,j$ of $I$ there exist an object $k$ and morphisms $k\to i$, $k\to j$.
\item For any two morphisms  $\alpha, \beta:i\to j$, there exists an object $k$ and a morphism $\gamma:k\to i$
such that $\alpha\circ \gamma=\beta\circ\gamma$.
\end{enumerate}
For each object $i$ of the small category $\cC$ one obtains a flat functor provided by the Yoneda embedding
$
h_i: \cC \longrightarrow \Se,~j\mapsto \Hom_\cC(i,j)$.
Here we take $\cC=\lbto$ and $F$ given by \eqref{archH1}. The filtering property of $F$ only involves finitely many elements of $\Hom_\frg(\frg(n),G)$ and hence by Lemma \ref{lemmain}, it follows using the filtering property of the functors $h_i$. The first part of Lemma \ref{lemmain} is used to prove the filtering property $(ii)$ while the last part is implemented to prove $(iii)$.\qed

%Subsection: The image of the flat functor by the module morphism
\subsection{The image of a flat functor $F:\lbto\longrightarrow \Se$ by the module morphism}\label{sectmod}

The functor $\Mod:\Arc\ltimes\N^\times\longrightarrow \nt$  associates to any morphism of archimedean sets the integer $k\in \nt$ involved in Definition \ref{last2} (\cf \cite{CCproj}). Its restriction  to the full subcategory whose objects are  archimedean sets of type $(\Z, \theta)$, where $\theta(x)=x+n+1$, defines a functor 
$\Mod:\lbt\longrightarrow \nt$. The category $\nt$ is isomorphic to its opposite in view of the commutativity of the multiplicative mono\"id  of positive integers.
The functor ${\rm Mod}:\lbt^{\rm op}\longrightarrow \N^\times$ determines a geometric morphism
of topoi. 
We recall once again (\!\!\cite{CCas}) that the category of points of the topos $\widehat{ \N^\times} $ is canonically equivalent to the category of totally ordered groups isomorphic to non-trivial subgroups of   $(\Q,\Q_+)$, and injective morphisms of ordered groups. This latter category is in turn equivalent  to the category  of algebraic extensions of the semifield  $\F=\Z_{\rm max}$ \ie of extensions 
$
\F\subset K\subset \Q_{\rm max}
$.
The morphisms are the {\em injective} morphisms of semifields.

This section is devoted to the description of the action of the geometric morphism $\Mod:\epi\longrightarrow \wnt$ on points, in terms of the associated flat functors. This process allows one to recover the extension $K$ of $\F=\Z_{\rm max}$ involved in Theorem \ref{main} from the datum of a flat functor $\lbt^{\rm op}\longrightarrow \Se$.

Given two small categories $\cC_j~(j=1,2)$, a functor $\phi:\cC_1\longrightarrow \cC_2$ determines a geometric morphism (also noted $\phi$) of topoi $\hat\cC_1\longrightarrow\hat\cC_2$ (\cf~\eg\cite{MM}, Chapter VII \S 2, Theorem 2). The inverse image $\phi^*$ sends an object of $\hat\cC_2$, \ie a contravariant functor $\cC_2\longrightarrow\Se$ to its composition with $\phi$ which determines a contravariant functor  $\cC_1\longrightarrow\Se$.
 The geometric morphism $\phi$ sends points of $\hat\cC_1$ to points of $\hat\cC_2$. In terms of the flat functors associated to points, the image by $\phi$ of a flat functor $F_1:\cC_1\longrightarrow \Se$ associated to a point $p_1: \Se\longrightarrow\hat\cC_1$ is the flat functor $F_2:\cC_2\longrightarrow \Se$ obtained by composing the Yoneda embedding $\cC_2\longrightarrow \hat \cC_2$ with $p_1^*\circ \phi^*$, where $p_1^*: \hat\cC_1\longrightarrow\Se$ is the inverse image functor with respect to $F_1$. Thus, for any object $Z$ of $\cC_2$ one obtains
\begin{equation}\label{yoneda1}
F_2(Z)=p_1^*(X), \   \   X:\cC_1^{\rm op}\longrightarrow \Se, \  \  X(c_1):=\Hom_{\cC_2}(\phi(c_1),Z).
\end{equation}
We apply this procedure to the functor ${\rm Mod}:\lbto\longrightarrow \N^\times$: \ie we take $\cC_1=\lbto$, $\cC_2=\nt$ and $\phi=\Mod$.
Let $F$ be a flat functor $F:\lbto\longrightarrow \Se$. The inverse image functor with respect to $F$ of a covariant functor $X:\lbt \longrightarrow \Se$ coincides  with the geometric realization $\vert X\vert_F$  and it 
is of the form
$$
   |X|_F=\left(\coprod_{n\ge 0} (F(\underline n)\times_\lbt  X(\underline n))\right)/\sim
$$
The image of the flat functor $F:\lbto\to \Se$ by the morphism $\mmod$ is thus the flat functor 
$H:\N^\times\longrightarrow \Se$ obtained as the geometric realization $\vert X\vert_F$ of the covariant functor
$$
X:\lbt \longrightarrow \Se,\qquad X(\underline n)=\Hom_{\nt}(\bullet, \mmod(\underline n))\cong \nt\qquad {\rm Obj}(\N^\times) = \{\bullet\}.
$$
The functor $X$ associates to any object of $\lbt$ the set $\nt$ and to a morphism  $\gamma$ of $\lbt$ its module $\mmod(\gamma)$ acting by multiplication on $\nt$. Hence we obtain 
\begin{equation}\label{H}
H(\bullet)=\left(\coprod_{n\ge 0} (F(\underline n)\times_\lbt  \nt)\right)/\sim
\end{equation} 
The equivalence relation is exploited as follows: for $(z,k)\in F(\underline n)\times  \nt$, one has $(z,k)\sim (F(\gamma)z,k)\in F(\underline 0)\times \nt$ for any $\gamma\in \Hom_\Lambda(\underline 0, \underline n)$, since $\Mod(\gamma)=1$. Moreover for $(z,k)\in F(\underline 0)\times  \nt$, one has $(z,k)\sim (F(\gamma)z,1)$ for 
$\gamma\in \Hom_\lbt(\underline 0, \underline 0)$ with $\mmod(\gamma)=k$. This shows that any element of $H(\bullet)$ is equivalent to an element of the form $(z,1)$ for some $z\in F(\underline 0)$. In particular one deduces that $H(\bullet)$ is a quotient of $F(\underline 0)$.
\begin{lem} \label{imagemod} Let  $F:\lbto\longrightarrow \Se$ be a flat functor and $(H,H_+)$ the corresponding point of $\widehat \nt$ through the morphism $\mmod$. Then there is a canonical, surjective $\nt$-equivariant map 
$$
\pi:F(\underline 0)\to H_+, \  \   \pi(F(\gamma)z)=\mmod(\gamma)\pi(z), \ \forall \gamma \in \Hom_\lbt(\underline 0, \underline 0).
$$
Moreover, the equivalence relation $x \sim x' \iff \pi(x)=\pi(x')$ is given by
\begin{equation}\label{equiv}
x \sim x' \iff \exists z\in F(\underline 1), \  F(\delta_0)z=x, \  \  F(\delta_1)z=x'
\end{equation}
where $\delta_j\in \Hom_\Delta(\underline 0, \underline 1),~(j=1,2)$ are the two face maps.
\end{lem}
\proof Let $\pi:F(\underline 0)\to H(\bullet)$, $\pi(z)=(z,1)$ as in \eqref{H}. As remarked above, the map $\pi$ is surjective. Let $\gamma \in \Hom_\lbt(\underline 0, \underline 0)$, then one derives easily: 
$
\pi(F(\gamma)z)=(F(\gamma)z,1)\sim (z,\mmod(\gamma))=\mmod(\gamma)\pi(z)
$. 
Moreover, $\pi(x)=\pi(x')$ if and only if there exist $n\in \N$, $u\in F(\underline n)$ and $\gamma, \gamma'\in \Hom_\lbt(\underline 0, \underline n)$ such that 
$
F(\gamma)u=x, \  \  F(\gamma')u=x', \  \mmod(\gamma)=\mmod(\gamma').$
Let $k=\mmod(\gamma)$, and $\Psi_k(\underline n)$  the archimedean set obtained from $\underline n=(\Z, \theta)$ by replacing $\theta$ with $\theta^k$ (\cf~Appendix~\ref{appcyclic}). Then one derives canonical factorizations involving the identity map $\Z\to\Z$ viewed as the element $\id_n^k\in \Hom_\lbt(\Psi_k(\underline n), \underline n)$ with $\mmod(\id_n^k)=k$
$$
\gamma=\id_n^k\circ \alpha_0, \  \  \gamma'=\id_n^k\circ \alpha_1, \  \   \alpha_0, \alpha_1\in \Hom_\Lambda(\underline 0, \Psi_k(\underline n)).
$$
 Let $\delta_j\in \Hom_\Delta(\underline 0, \underline 1)$, $j=0,1$, be the two face maps. One can then find  $\alpha\in \Hom_\Lambda(\underline 1, \Psi_k(\underline n))$ such that $\alpha_j=\alpha\circ \delta_j$. Thus it follows that
$
\gamma=\id_n^k\circ \alpha\circ \delta_0, \    \gamma'=\id_n^k\circ \alpha\circ \delta_1.
$
and one also gets that 
$
x=F(\gamma)z=F( \delta_0)F(\alpha)F(\id_n^k)u=F( \delta_0)z,   \ z=F(\alpha)F(\id_n^k)u, \ \  x'=F(\delta_1)z,
$
which proves \eqref{equiv}.\qed

Let $\iota: \nt\longrightarrow \lbto$ be the functor which associates to the unique object $\bullet$ of the small category $\nt$ the object $\underline 0$ of $\lbto$ and to a positive integer $k$ the unique element $\iota(k)\in \Hom_\lbt(\underline 0, \underline 0)$ such that $\Mod(\iota(k))=k$.
There is another natural equivalence relation on $F(\underline 0)$ deduced from the action of $\nt$ and explicitly given by 
\begin{equation}\label{Q}
x\sim' y \iff \exists n, m \in \nt, \ \ F(\iota(n))x=F(\iota(m))y.
\end{equation}
We let $\pi':  F(\underline 0)\to Q$ be the map to the  quotient $Q:=F(\underline 0)/\sim'$.
\begin{lem} Let  $F:\lbto\longrightarrow \Se$ be a flat functor. Let $(H,H_+)$ be the corresponding point of $\widehat \nt$ through the morphism $\mmod$. Then, with $\pi':  F(\underline 0)\to Q$ defined just above, the map 
$$
(\pi,\pi'):F(\underline 0)\to H_+\times Q
$$
is a bijection of sets.
\end{lem}
\proof We prove first that $(\pi,\pi')$ is injective. Let $x,y\in F(\underline 0)$ be such that $x\sim' y$. Let $n, m \in \nt$ with  $F(\iota(n))x=F(\iota(m))y$. Then the $\nt$-equivariance of $\pi$ shows that in the ordered group $H$ one has $n\pi(x)=m\pi(y)$. Thus the equality $\pi(x)=\pi(y)$ implies that $n=m$ since $H$ has no torsion and both $\pi(x)$ and $\pi(y)$ are non-zero. Hence we derive $F(\iota(n))x=F(\iota(n))y$. So far we have used $F:\lbto\longrightarrow \Se$ as a contravariant functor from $\lbt$ to $\Se$. In this proof 
we shall also use $F$ directly as a covariant functor and for clarity of exposition we shall denote  $F^{\rm cov}(\alpha)=F(\alpha^*)$ for any morphism $\alpha\in \lbto$ using the canonical anti-isomorphism $\lbt\to \lbto$. The flatness of $F$ implies that there exists an object $(\underline m, t)$ of $\int_\lbto F$, with $t\in F(\underline m)$, and two elements $u,v\in \Hom_\lbto(\underline m, \underline 0) $ such that $F^{\rm cov}(u)t=x$ and $F^{\rm cov}(v)t=y$.  One thus obtains $F^{\rm cov}(\iota(n)^*\circ u)t=F^{\rm cov}(\iota(n)^*\circ v)t$. The flatness of $F$ provides an object $(\underline \ell,s)$ of $\int_\lbto F$, with $s\in F(\underline \ell)$, and a morphism  in $\int_\lbto F$, given by a $w\in \Hom_\lbto(\underline \ell, \underline m) $ such that $\iota(n)^*\circ u\circ w=\iota(n)^*\circ v\circ w$. The fact that $w$ is a morphism  in $\int_\lbto F$ connecting the object $(\underline \ell,s)$ to $(\underline m, t)$ means that $F^{\rm cov}(w)s=t$. The equality $\iota(n)^*\circ u\circ w=\iota(n)^*\circ v\circ w$ in $\lbto$ implies that the two elements $u\circ w$, $ v\circ w$ of $\Hom_\lbto(\underline \ell,\underline 0)$ are equal.  Indeed, let us show that for $\gamma, \gamma'\in \Hom_\lbto(\underline \ell,\underline 0)$, and $n\in \nt$
\begin{equation}\label{injectell}
\iota(n)^*\circ\gamma=\iota(n)^*\circ\gamma'\implies \gamma=\gamma'.
\end{equation}
To prove this we pass to the opposite category using the identification $\Hom_\lbto(\underline \ell,\underline 0)\cong  \Hom_\lbt(\underline 0,\underline \ell)$, thus we get elements $\alpha,\alpha'\in \Hom_\frg(\frg(0),\frg(\ell))$ such that $\alpha\circ \iota(n)=\alpha'\circ \iota(n)$. 
The \totg $\frg(0)$ has a single object $o$ and an element  $\alpha\in\Hom_\frg(\frg(0),\frg(\ell))$ is uniquely specified by the image $\alpha(o)$ and by its module $\Mod(\alpha)$. Replacing $\alpha$ by  $\alpha\circ \iota(n)$ does not alter the image $\alpha(o)$ while it replaces the module by $n\,\Mod(\alpha)$. Thus \eqref{injectell} follows. This shows that  $u\circ w=v\circ  w$ and $x=F^{\rm cov}(u\circ w)s=F^{\rm cov}(v\circ w) s=y$ which proves that $(\pi,\pi')$ is injective.

Next we show that $(\pi,\pi')$ is also surjective. It is enough to prove that given $x,y\in F(\underline 0)$ there exists  $z\in F(\underline 0)$ such that
\begin{equation}\label{Qs}
x\sim' z, \   \  y\sim z.
\end{equation}
Again from the flatness of $F$ one derives that there exists an object $(\underline m,t)$ of $\int_\lbto F$, with $t\in F(\underline m)$, and two elements $u,v\in \Hom_\lbto(\underline m, \underline 0) $ such that $F^{\rm cov}(u)t=x$ and $F^{\rm cov}(v)t=y$. Let then $u^*, v^*\in \Hom_\lbt(\underline 0, \underline m)$ be the corresponding morphisms in the dual category $\lbt$. We view them as elements of $\Hom_\cG(\frg(0),\frg(m))$. Then there exists a unique $w^*\in \Hom_\cG(\frg(0),\frg(m))$ which has the same range object of $\frg(m)$ as $u^*$ and the same module as $v^*$.
Since $w^*$ and $v^*$ have the same module $k$ one can find, as in the proof of Lemma \ref{imagemod},  $\alpha\in \Hom_\Lambda(\underline 1, \Psi_k(\underline m))$ such that 
$
w^*=\id_m^k\circ \alpha\circ \delta_0, \   v^*=\id_m^k\circ \alpha\circ \delta_1
$.

Let $z=F^{\rm cov}(w)t$, one has
$
z=F(\delta_0)F(\id_m^k\circ \alpha)t, \    y=F(\delta_1)F(\id_m^k\circ \alpha)t
$
which shows that $y\sim z$. Since $w^*\in \Hom_\cG(\frg(0),\frg(m))$ has the same range object of $\frg(m)$ as $u^*$, there exist $k,k'\in \nt$ such that $u^*\circ \iota(k)=w^*\circ \iota(k')$. This implies that 
$
F( \iota(k))x=F^{\rm cov}(\iota(k)^*)F^{\rm cov}(u)t=F^{\rm cov}(\iota(k')^*)F^{\rm cov}(w)t=F( \iota(k'))z 
$,
so that $x\sim' z$ as required.\qed

\begin{cor} Let  $F:\lbto\longrightarrow \Se$ be a flat functor. Then, for any $x\in F(\underline 0)$ the map $\pi: F(\underline 0)\to H_+$ induces a bijection between the  commensurability class of $x$ under the action of $\nt$ and the positive part $H_+$ of the ordered group image of the point $F$ by the geometric morphism $\Mod:\epi\longrightarrow \wnt$.
\end{cor}

\begin{figure}
\begin{center}
\includegraphics[scale=0.4]{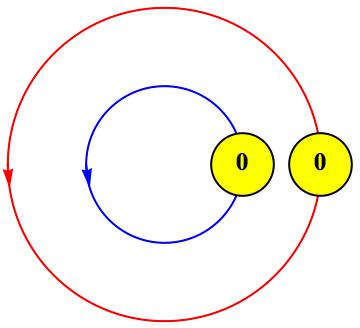}
\end{center}
\caption{The morphism $\frg(0)\to \frg(n)$ is the unique element $j_n\in\Hom_\dop(0^*,n^*)\subset \Hom_\lbt(\underline 0,\underline n)$. It gives a map $F(j_n):F(\underline n)\to F(\underline 0)$ whose value fixes an object and the length of the interval.\label{clock2} }
\end{figure}

%Subsection: The interval associated to a pair...
\subsection[The interval associated to a pair $F:\lbto\longrightarrow \Se$, $x\in F(0)$]{The interval associated to a pair $F:\lbto\longrightarrow \Se$, $x\in F(\underline 0)$}\label{intervalsubsec}

In this section we describe a process which associates an interval $F_x$ to a pair made by a flat functor $F:\lbto\longrightarrow \Se$ and an element $x\in F(\underline 0)$. For each integer $n\geq 0$ we denote by $j_n\in\Hom_\dop(0^*,n^*)\subset \Hom_\lbt(\underline 0,\underline n)$ the unique element of $\Hom_\dop(0^*,n^*)$. One has $j_n(x)=(n+1)x$, $\forall x\in \Z$. We refer to Appendix \ref{appbary} for the notations on the category $\dop$. 
We identify the opposite $\dop$ of the simplicial category $\Delta$ with the subcategory of the epicyclic category $\lbt$ with the same objects and morphisms specified by the identification
\begin{equation}\label{mordop}
\Hom_\dop(n^*,m^*)\simeq \{\gamma \in \Hom_\lbt(\underline n,\underline m)\mid \gamma\circ j_n=j_m \}.
\end{equation}
Passing to the dual categories, we get $j_n^*\in \Hom_\lbto(\underline n,\underline 0)$ and these elements define objects of the slice category
$\lbto/\underline 0$.  More precisely one obtains

\begin{lem}\label{lemslice} The simplicial category $\Delta$ is canonically isomorphic to the full subcategory of the slice category
$\lbto/\underline 0$ with objects the $j_n^*$.
\end{lem}
 \proof We define  a functor $\kappa:\Delta\longrightarrow\lbto/\underline 0$ sending the object $[n]=\{0,\ldots,n\}$ to $j_n^*$. To define the action of $\kappa$
on morphisms we first describe explicitly the identification \eqref{mordop}. Given $\phi\in \Hom_\dop(n^*,m^*)$ \ie a non-decreasing map preserving the end points from the interval $n^*:=\{0, \ldots , n+1\}$ to the interval $m^*$, one extends $\phi$ uniquely to the following non-decreasing map 
$$
\tilde \phi: \Z  \to \Z,\  \   \tilde \phi(x)=\phi(x) \qqq x\in \{0, \ldots , n+1\}, \ \tilde \phi(x+a(n+1))=\tilde\phi(x)+a(m+1)\qqq a\in \Z.
$$
One finds $\tilde \phi \in \Hom_\lbt(\underline n,\underline m)$ with $\Mod(\tilde \phi)=1$. The map $\phi\mapsto \tilde \phi$ defines a faithful, covariant functor $\dop\longrightarrow\lbt$ which is bijective on objects and whose range is characterized by \eqref{mordop}. To check this statement, note that the equation $\gamma\circ j_n=j_m$ implies that $\Mod(\gamma)=1$. Moreover any non-decreasing map $\psi: \Z  \to \Z$ such that $\psi(0)=0$ and $\psi(x+a(n+1))=\psi(x)+a(m+1),~\forall a\in \Z$ is of the form $\tilde \phi$ for a unique $\phi\in \Hom_\dop(n^*,m^*)$.
Then we define $\kappa$ on morphisms as follows 
$$
\kappa(\gamma):=\left(\tilde\gamma^*\right)^*\in \Hom_{\lbto/\underline 0}(j_n^*,j_m^*)\qqq \gamma \in \Hom_\Delta([n],[m]).
$$\qed

\begin{lem}\label{leminterval}
Let  $F: \lbt^{\rm op}\longrightarrow \Se$ be a flat functor and $x\in F(\underline 0)$. The following equality defines a flat functor  $\Delta\longrightarrow \Se$ (\ie interval) which can be equally interpreted as a contravariant  functor 
\begin{equation}\label{mordop1}
F_x:\dop \longrightarrow \Se,\qquad F_x(n^*)=\{z\in F(\underline n)\mid F(j_n)z=x\}.
\end{equation}
\end{lem}
\proof Let $X$ be a totally ordered set on which an ordered group $H$ acts fulfilling the conditions of Lemma \ref{lemmain}.
Assume first that the functor  $F$ is defined by
$$
F(\underline n)=\Hom_\frg(\frg(n),G), \  \ G=G(X,H).
$$
Then, by applying \eqref{functn}, $F(\underline n)$ is the quotient by the diagonal action of $H$ on the subset  of $X^{n+2}$
$$\{(x_0,\ldots , x_{n+1})\in X^{n+2}\mid x_j\leq x_{j+1}, \forall j\leq n,~ x_{n+1}= x_0+ H_+\}.
$$
Thus one derives that 
$
F(\underline 0)=X/H\times H_+.
$
Let  $x\in F(\underline 0)$, $x=(\tilde x_0,h)\in X/H\times H_+$. After choosing a lift $x_0\in X$ of $\tilde x_0\in X/H$,   \eqref{mordop1} provides the equality
\begin{equation}\label{interval}
F_x(n^*)=\{(x_0,\ldots , x_{n+1})\in X^{n+2}\mid x_j\leq x_{j+1}, \forall j\leq n, ~x_{n+1}= x_0+h\}
\end{equation}
which can be equivalently described as $\Hom(n^*,I)$, where $I$ is the interval $[x_0,x_0+h]\subset X$. Let $m\in \N$ and $G=\frg(m)$, 
 $F(-)=\Hom_\lbt(-, \frg(m))$.  Then $x\in F(\underline 0) $ is determined by an object of $\frg(m)$ and an integer $k$, and it follows by the above discussion and from the arguments developed in Appendix~\ref{appbary} that the associated  interval is $F_x=\sdd_k(m^*)$. \newline
The general case is deduced from the above one by writing the  flat functor $F$ as a filtering colimit of functors of the form $\Hom_\lbt(-, \frg(m))$.
More precisely, let $F: \lbt^{\rm op}\longrightarrow \Se$ be a flat functor and $x\in F(\underline 0)$. We prove, for instance, the  filtering property of $F_x$. An object $\alpha$ of the category $\int_\Delta F_x$ is of the form 
$
\alpha=(n,z), \  z\in F(\underline n), \   F(j_n)z=x.
$
Given two such objects $\alpha_j=(n_j,z_j),~j=1,2$ one has $F(j_{n_1})z_1=F(j_{n_2})z_2$ and since $F$ is a filtering colimit of functors of the form $\Hom_\lbt(-, \frg(m))$ one can find an integer $m\in \N$ and elements $y_j\in \Hom_\lbt(\frg(n_j), \frg(m))$ such that $y_j$ represents $z_j$ in the colimit and that $y_1\circ j_{n_1}=y_2\circ j_{n_2}$. This last equality can be realized as a consequence of $F(j_{n_1})z_1=F(j_{n_2})z_2$ using the filtering colimit \ie the definition of equality in the limit. The existence of an object $\alpha=(n,z)$ and morphisms $u_j$ in $\int_\Delta F_x$ to $\alpha_j=(n_j,z_j)$ then follows from the filtering property of the category 
$\int_\Delta F'_y$ where $F'$ is the flat functor $\Hom_\lbt(-, \frg(m))$ and $y=y_1\circ j_{n_1}=y_2\circ j_{n_2}$.\qed

%Subsection: The oriented groupoid associated to a pair...
\subsection[The \totg associated  to a pair $F:\lbto\longrightarrow \Se$, $x\in F(0)$]{The \totg associated  to a pair $F:\lbto\longrightarrow \Se$, $x\in F(\underline 0)$}\label{recsect}

Let $F:\lbto\longrightarrow \Se$ be a flat functor.
Lemma \ref{leminterval} shows how to associate an interval $F_x$ to an element $x\in F(\underline 0)$. Next, we state two lemmas which will be used to reconstruct the totally ordered set $X$ with an action of $(H,H_+)$, where $(H,H_+)$ is  the ordered group corresponding to $F$ by means of the geometric morphism $\Mod: \epi \longrightarrow \widehat{\N^\times}$, as a point of $\wnt$: \cf ~Lemma \ref{imagemod}. We use the notations of Appendix~\ref{appbary}.
 
\begin{lem}\label{lemsubint}
Let $F: \lbt^{\rm op}\longrightarrow \Se$ be a flat functor and $x\in F(\underline 0)$. The following map defines an isomorphism of intervals for any  $k\in \N^\times$\begin{equation}\label{explicitis0}
\omega_k:\sdd_k(F_x)\to F_{F(\iota(k))x},\quad \omega_k(z,\rho)=F(\id_m^k\circ \rho)z\in F_{F(\iota(k))x}(n^*)
\end{equation}
$\forall (z,\rho)\in F_x(m^*) \times \Hom_\dop(n^*,\sdd_k(m^*))$.
\end{lem}
\proof Using the same method as in the proof of Lemma \ref{leminterval}, it is enough to check that the map is well defined and that one gets an isomorphism in the special case where the functor $F$ is of the form $\Hom_\lbt(\bullet, \frg(m))$. We consider the subcategory $\dop\ltimes \N^\times\subset \lbt$ also reviewed  in Appendix~\ref{sectcrossprod}. 

\begin{figure}
\begin{center}
\includegraphics[scale=0.6]{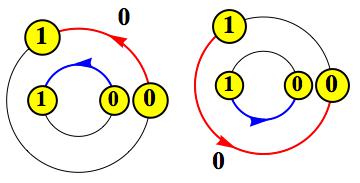}
\end{center}
\caption{The morphism $\eta:\frg(1)\to \frg(2)$ is in $\Hom_\dop(1^*,2^*)\subset \Hom_\lbt(\underline 1,\underline 2)$. It determines a  map $F(\eta):F(\underline 2)\to F(\underline 1)$ such that $F(\eta)F_x(2^*)\subset F_x(1^*)$, since $\eta \circ j_1=j_2$.\label{clock3} }
\end{figure}

In Appendix~\ref{sectexplicit}, Lemma \ref{lemdescsub}, it is shown that the action of the barycentric subdivision functor $\sd_k$ on points of the topos of simplicial sets $\hat \Delta$ is described in terms of the associated flat functors by replacing the given flat functor $G:\Delta\longrightarrow \Se$ with the flat functor 
$\sdd_k(G):\Delta\longrightarrow \Se$. This latter associates to the object $[n]$ of $\Delta$ the quotient set
\begin{equation}\label{explicit0}
  \sdd_k(G)([n])=\left(\coprod_{m\ge 0} (G(m)\times \Hom_\dop(n^*,\sdd_k(m^*)))\right)/\sim
\end{equation}
where the equivalence relation $\sim$ is given, for $f\in \Hom_\dop(m^*,r^*)=\Hom_\Delta([r],[m])$, $\beta\in \Hom_\dop(n^*,\sdd_k(m^*))$ and $y\in G(r)$ as follows
$$G(m)\times \Hom_\dop(n^*,\sdd_k(m^*))\ni
(G(f)y,\beta)\sim(y,\sdd_k(f)\circ \beta)\in G(r)\times \Hom_\dop(n^*,\sdd_k(r^*)).
$$
We prove that the map \eqref{explicitis0} is well defined.  First, one has $\id_m^k\circ \rho\in \Hom_\lbt(\underline n,\underline m)$ and since $\id_m^k\circ \rho\circ j_n=j_m\circ \iota(k)$ one gets
$$
F(\id_m^k\circ \rho)z \in F(\underline n), \quad  F(j_n)F(\id_m^k\circ \rho)z=F(\id_m^k\circ \rho\circ j_n)z=F(\iota(k))x.
$$
We let $G(m)=F_x(m^*)$. We  need to check that two equivalent elements in \eqref{explicit0}
have the same image under $\omega_k$. Let $f\in \Hom_\dop(m^*,r^*)=\Hom_\Delta([r],[m])$, $\beta\in \Hom_\dop(n^*,\sdd_k(m^*))$ and $y\in G(r)=F_x(r^*)$. One then has, since $f\circ\id_m^k=\id_r^k\circ\sdd_k(f)$
$$
\omega_k(G(f)y,\beta)=F(\id_m^k\circ \beta)G(f)y=F(\id_m^k\circ \beta)F(f)y=F(f\circ\id_m^k\circ \beta)y=F(\id_r^k\circ\sdd_k(f)\circ  \beta)y.
$$
Thus one derives 
$
\omega_k(G(f)y,\beta)=F(\id_r^k\circ\sdd_k(f)\circ  \beta)y=\omega_k(y,\sdd_k(f)\circ \beta)
$.
This shows that the map $\omega_k$ is well defined. One easily checks that it is an isomorphism in the special case when  the functor $F$ is of the form $\Hom_\lbt(\bullet, \frg(m))$. Then,  one obtains the general case working as in the proof of Lemma \ref{leminterval}.\qed \endproof

Given a flat functor $H: \Delta\longrightarrow\Se$,   the set $I:=H([1])$ is endowed with the total order relation defined by 
$$
y \leq y' \iff \exists z\in H([2]) \mid y=H(\sigma)(z), \, y'=H(\sigma')(z)
$$
where $\sigma, \sigma'\in \Hom_\Delta([2],[1])$ are the two surjections.  Endowed with this order relation $I$ is an interval and one recovers the flat functor $H$ as $H([n])=\Hom_\dop(n^*,I)$. We shall use this fact to organize the intervals $I_x=F_x([1])$ associated to the pair of a flat functor $F:\lbto\longrightarrow \Se$ and $x\in F(\underline 0)$.\newline
Our final goal is in fact to show how to reconstruct an oriented groupo\"id from a flat functor $\lbt^{\rm op}\longrightarrow \Se$. It will be enough to check that we can reconstruct the oriented groupo\"id $G(X,H)$ for flat functors of the form $F(\underline m)=\Hom_\frg(\frg(m),G(X,H))$ as in Proposition \ref{propmain}. 

The $k$ elements of $\Sigma_1^k$ (\cf Appendix A) determine the following elements  $\psi_{k,j}\in \Hom_{\dop\ltimes \N^\times}(1^*,1^*)\subset \Hom_\frg(\frg(1),\frg(1))$, for $0\leq j<k$,
$$
\psi_{k,j}:\ps_1^k\circ \alpha_{k,j}, \  \alpha_{k,j}\in \Hom_\dop(1^*,\sdd_k(1^*)), \quad \alpha_{k,j}(1)=2j+1.
$$
and we view them as the elements of $\Hom_\frg(\frg(1),\frg(1))$ associated to the endomorphisms $\phi_{k,j}$ (as in Definition \ref{last2})  of the archimedean set $(\Z,u\mapsto u+2)$ given  by
\begin{equation}\label{psikj}
\phi_{k,j}: \Z\to \Z, \ \  \phi_{k,j}(2\ell):=2k\ell, \ \ \phi_{k,j}(2\ell+1):=2j+1+2k\ell
\end{equation}
Note that $\mmod(\psi_{k,j})=k$. One lets $\psi_k:=\psi_{k,1}$.

Note that any $\psi \in \Hom_\lbt(\frg(1),\frg(1))=\Hom_\frg(\frg(1),\frg(1))$ such that $\psi\circ j_1=j_1\circ \iota(k)$ defines a map 
\begin{equation}\label{explicit3}
F_\psi:I_x\to I_{F(\iota(k))x},\quad F_\psi(z)=F(\psi)z, \ \forall z \in F(\underline 1), \ F(j_1)z=x.
\end{equation}
Indeed one has $F(\psi)z\in F(\underline 1)$ and $F(j_1)F(\psi)z=F(\iota(k)F(j_1)z=F(\iota(k))x$.  Thus $F(\psi)z\in I_{F(\iota(k))x}$. This fact is used in the statement of the next proposition.

\begin{prop}\label{propreconstruct} Let $F$ be a  flat functor of the form $F(\underline m)=\Hom_\frg(\frg(m),G(X,H))$ as in Proposition \ref{propmain}. Let $x\in F(\underline 0)$, $x=(\tilde x_0,h)\in X/H\times H_+$ and $x_0\in X$ be a lift of $\tilde x_0$. 

$(i)$~The following equality determines the ordered set $Y=\{y\in X\mid y\geq x_0\}$:
\begin{equation}\label{defnJ}
J_x:=\varinjlim \left(I_{F(\iota(k))x}, F_{\psi_{\ell/k}}\right)
\end{equation}
where the colimit is taken using the maps  $F_{\psi_k}:I_x\to I_{F(\iota(k))x}$ and the set of indices is ordered by divisibility.\vspace{0.1in}

$(ii)$~The graph of the action of $h$ on $Y$ is given as the union in the inductive limit of the subsets 
\begin{equation}\label{defnJbis}
\bigcup_{0\leq j <k-1} \{(F_{\psi_{k,j}}(z),F_{\psi_{k,j+1}}(z))\mid z\in I_x\}
\end{equation}
\vspace{0.1in}

$(iii)$~ The \totg $G(X,H)$ is isomorphic to the \totg of ordered pairs of elements of $J_x$ modulo 
the action of $H$.
\end{prop}
\proof By \eqref{functn}, \eqref{functnbis}, a morphism $\phi\in \Hom_\frg(\frg(m),G(X,H))$ is determined by a non-decreasing map $f:\Z\to X$ and an element $h\in H_+$, $h\neq 0$  such that $f(x+m+1)=f(x)+h$ $\forall x\in \Z$. Moreover, by \eqref{functn}, the maps $f$ and $g$ define the same morphism if and only if there exists an $a\in H$ such that $g=f+a$. In particular the map $f:\Z\to X$ associated to an element $z$ of the interval $I_x=[x_0,x_0+h]$ of \eqref{interval} is given by $f(2n)=x_0+nh$, $f(2n+1)=z+nh$ $\forall n\in \Z$. With these notations let us determine the map $F_{\psi_{k,j}}$. It corresponds to the composition $f\circ \phi_{k,j}$ and is hence given by
$$
f\circ \phi_{k,j}(0)=f(0)=x_0, \  \  f\circ \phi_{k,j}(1)=f(2j+1)=z+jh, 
$$
while for any $\ell\in \Z$ one has $f\circ \phi_{k,j}(u+2\ell)=f\circ \phi_{k,j}(u)+k\ell h$. Thus the map $F_{\psi_{k,j}}$ sends $z\in [x_0,x_0+h]$ to $z+jh\in [x_0,x_0+kh]$. In particular the $F_{\psi_{\ell/k}}$ give the canonical inclusion of  intervals $[x_0,x_0+kh]\subset [x_0,x_0+\ell h]$ for $k\vert \ell$ and the inductive limit \eqref{defnJ} is $Y=\{y\in X\mid y\geq x_0\}$. Moreover, in the limit, the graph of the action of $h$ is given by the pairs \eqref{defnJbis}. Finally one reconstructs in this way the action of $H_+$ on the  colimit set $J_x$ as well as the \totg of ordered pairs of elements of $J_x$ modulo 
the equivalence relation generated by the action of the semigroup $H_+$.
 \qed.

\begin{figure}
\begin{center}
\includegraphics[scale=0.7]{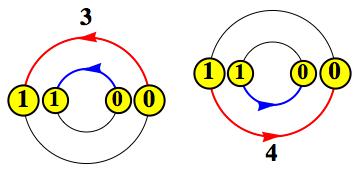}
\end{center}
\caption{The morphism $\psi:\frg(1)\to \frg(1)$ has module $8$, it gives a map $F_\psi:F_x\to F_{F(\iota(8))x}\cong\sdd_8(F_x)$.\label{clock} }
\end{figure}

%Proof of Theorem
\subsection{Proof of Theorem \ref{main}}

In Proposition \ref{propmain} we have shown  that a pair $(K,E)$ of an algebraic extension $K$ of $\F=\Z_{\rm max}$ and   an archimedean \mod~$E$ over $K$, defines a flat functor $\lbto\longrightarrow\Se$ as in \eqref{archH1}. Conversely, in \S \ref{sectmod} and  Proposition \ref{propreconstruct}, we have proven that a flat functor $F: \lbto\longrightarrow\Se$ determines the data given by a totally ordered group $(H,H_+)$ of rank one, a totally ordered set $J$ and an action of $H$ on $J$. One easily checks that the pair $(J,H)$ fulfills the conditions of Proposition \ref{propmain}  and that
\begin{equation*}%\label{archH3}
F(\underline n)=\Hom_\frg(\frg(n),G), \  \ G=G(J,H).
\end{equation*}
The proof of the Theorem~\ref{main} is then complete.

%Subsection: Exotic points
\subsection{Exotic points}

In general, a colimit involves the process of taking a  quotient and there is no guarantee ``a priori'' that the quotient is a non-singular space, even at the set-theoretic level. Next, we exhibit certain points of the epicyclic topos whose associated flat functors show a singular behavior.

Consider $K=\Q_{\rm max}$  and the \mod~$E=\R_{\rm max}$. Then the associated functor $F$ is such that $F(\underline n)$ is the quotient of the following subset of $\R^{n+2}$ by the diagonal action of $\Q$
$$\{(x_0,\ldots , x_{n+1})\in \R^{n+2}\mid x_j\leq x_{j+1}, \forall j\leq n,~ x_{n+1}- x_0\in\Q_+\}.
$$
One gets for instance
$$
F(\underline 0)=\R/\Q\times \Q_+.
$$
The units of the associated groupo\"id  $G=G(E,K)$ form the quotient set $\R/\Q$ which is singular. This result supports the point of view that the  description of the points of the epicyclic topos in terms of the category  $\cP$ as in Theorem \ref{main} is more appropriate than the  description given in terms of \totgs since the latter interpretation involves singular quotients.

\subsection{Relation between  $\widehat{\Lambda^{\rm op}}$ and $\epi$}

By a result of \cite{Moerdjik} the points of the topos $\widehat{\Lambda^{\rm op}}$ dual to the (opposite of) the cyclic category form the category of abstract circles and in turns by \cite{topos} this latter is equivalent to the category $\Arc$ of archimedean sets with morphisms of module $1$. The inclusion of categories $\Lambda^{\rm op}\subset \lbto$ (dual to $\Lambda\subset \lbt$) induces a geometric morphism  of topoi: in this section we determine its behavior on the points. We start by stating a technical lemma which will be applied in the proof of Proposition~\ref{effectinclus}

\begin{lem}\label{lemdecfrg} Let $(X,\theta)$ be an archimedean set and $\phi\in \Hom_\frg(\frg(u), G(X,\theta))$, $u\geq 0$. Denote by $r+1$  the cardinality of $\phi(\frg(u)^{(0)}) \subset G(X,\theta)^{(0)}$. Then there exists $g\in \Hom_\frg(\frg(u), \frg(r))$
and $f\in \Hom_\Arc(\underline r, (X,\theta))$ such that $\phi=G(f)\circ g$. This factorization is unique up to the cyclic group $\Aut_\Arc(\underline r)$, \ie the replacement 
$$
f\mapsto f\circ \sigma^{-1}, \  \  g\mapsto G(\sigma)\circ g\qqq \sigma\in \Aut_\Arc(\underline r).
$$
\end{lem}
\proof
 Let $\phi\in \Hom_\frg(\frg(u), G(X,\theta))$ and $k=\Mod(\phi)$. Then by 
\eqref{functn}, $\phi$ is determined by $u+2$ elements $x_j$, $0\leq j\leq u+1$ of $X$ such that $x_j\leq x_{j+1}$ $\forall j\leq u$ and that $x_{u+1}=\theta^k(x_0)$. One has a canonical isomorphism $[x_0,\theta^k(x_0)]\simeq \sdd_k([x_0,\theta(x_0)]$. Let $\chi\in \Hom_\dop(u^*,\sdd_k([x_0,\theta(x_0)])$ be  the morphism of intervals  defined by  $\chi(j)=x_j$, $\forall j, 0\leq j \leq u+1$. Let $Y=\phi(\frg(u)^{(0)}) \subset G(X,\theta)^{(0)}= X/\theta$. Then the rank  of the morphism $\chi$ is the cardinality of $Y\setminus \{\tilde x_0\}$ where $\tilde x_0\in  X/\theta$ is the class of $x_0$ and it is thus equal to $r$. By Proposition \ref{propapp1} there exists  a unique decomposition 
$$
\chi=\sdd_k(\alpha)\circ \beta, \ \ \beta\in \Hom_\dop(u^*,\sdd_k(r^*)), \ \alpha\in \Hom_\dop(r^*,[x_0,\theta(x_0)]).
$$
 The map $\alpha\in \Hom_\dop(r^*,[x_0,\theta(x_0)])$ is uniquely determined by its restriction to $\{1, \ldots, r\}$ which is the unique increasing injection of $\{1, \ldots, r\}$ in $(x_0,\theta(x_0))$ whose range  gives  $Y\setminus \{\tilde x_0\}\subset X/\theta$. 
Let $f\in \Hom_\Arc(\underline r,(X,\theta))$ be defined by   
$$
f(i+\ell(r+1))=\theta^\ell(\alpha(i))\qqq i, 0\leq i\leq r, ~\forall \ell \in\Z
$$ 
Moreover $\beta\in \Hom_\dop(u^*,\sdd_k(r^*))$ determines the element $g=\id_r^k \circ \beta\in \Hom_\frg(\frg(u), \frg(r))$. One has by construction $\Mod(g)=k$. The equality $\chi=\sdd_k(\alpha)\circ \beta$ shows that 
$\phi=G(f)\circ g$. The uniqueness of the decomposition follows from the uniqueness of the archimedean set $(X',\theta)$ obtained as the inverse image of $Y$ under the canonical projection $X\to X/\theta$. This gives the uniqueness of $f$ up to cyclic permutations and since $f$ is injective the uniqueness of $g$ follows. \qed

\begin{prop}\label{effectinclus} The geometric morphism $j: \widehat{\Lambda^{\rm op}}\longrightarrow \epi$ associated to the inclusion of categories $\Lambda^{\rm op}\subset \lbto$ induces on the   points of the  topoi the functor producing the inclusion of $\F$-\mods~among  \mods~over algebraic extensions of $\F$.
\end{prop}
\proof A point $p$ of $\widehat{\Lambda^{\rm op}}$ is given by a flat functor of the form
$$
P:\Lambda^{\rm op}\longrightarrow\Se, \  \   \underline n \mapsto P(\underline n)=  \Hom_\Arc(\underline n,(X,\theta))
$$
where  $(X, \theta)$ is  an archimedean set and, as above, $\underline n=(\Z,t^{n+1})$ is the archimedean set given by the translation $x\mapsto x+n+1$ on $\Z$.  
By \eqref{yoneda1}, $j(p)$ is given by  the flat functor $F:\lbto\longrightarrow \Se$
\begin{equation}\label{yoneda2}
F(Z)=p^*(X), \   \   X:\Lambda\longrightarrow \Se, \  \  X(\underline n):=\Hom_{\lbto}(\frg(n),Z)=\Hom_{\frg}(Z,\frg(n)).
\end{equation}
The inverse image $p^*(X)$ is defined as 
$$
p^*(X)= \coprod P(\underline n)\times X(\underline n)/\sim\, = \coprod \Hom_\Arc(\underline n,(X,\theta))\times \Hom_{\frg}(Z,\frg(n))/\sim
$$
where the equivalence relation $\sim$ is the simplification with respect to $\Lambda$, \ie one has 
$$
(f\circ h,g)\sim (f,G(h)\circ g)\qqq h\in \Hom_\Arc(\underline n,\underline m),~ f\in P(\underline m),~ g\in X(\underline n).
$$
 It follows that the map $\gamma:p^*(X)\to \Hom_\frg(Z, G(X,\theta))$ which sends $(f,g)\in P(\underline n)\times X(\underline n)$ to $\gamma(f,g)=G(f)\circ g\in \Hom_\frg(Z, G(X,\theta))$ is well defined. Next, we show that when $Z=\frg(u)$ ($u\geq 0$) this map is bijective. \newline
First, we prove that $\gamma$  is surjective. Let $\phi\in \Hom_\frg(\frg(u), G(X,\theta))$. By Lemma \ref{lemdecfrg}, there exists $r\geq 0$, $g\in \Hom_\frg(\frg(u), \frg(r))$
and $f\in \Hom_\Arc(\underline r, (X,\theta))$ such that $\phi=G(f)\circ g$. Thus $\phi=\gamma(f,g)$ is in the range of $\gamma$.\newline
We prove that $\gamma$ is injective.  It is enough to show that for any $n\geq 0$ and any pair
$
(f',g')\in \Hom_\Arc(\underline n,(X,\theta))\times \Hom_{\frg}(\frg(u),\frg(n))
$ 
one has $(f',g')\sim (f,g)$, where the pair $(f,g)$ is obtained from Lemma \ref{lemdecfrg} applied to  $\phi=f'\circ g'\in \Hom_\frg(\frg(u), G(X,\theta))$. The element $g'\in \Hom_{\frg}(\frg(u),\frg(n))$ satisfies $\Mod(g)=\Mod(f'\circ g')=k$. One has $\frg(n)=G(X_0,\theta_0)$ where $X_0=\Z$ and $\theta_0(x)=x+n+1$.  Thus by Lemma \ref{lemdecfrg} applied to $g'\in \Hom_\frg(\frg(u), G(X_0,\theta_0))$ one obtains a decomposition  
 $$g'=G(f_0)\circ g_0, \  \  f_0\in \Hom_\Arc(\underline r_0,(X_0,\theta_0)), \ \  g_0 \in \Hom_\frg(\frg(u), \frg(r_0))$$
One has $f_0\in \Hom_\Arc(\underline r_0,\underline n)$ and thus $(f',g')\sim (f'\circ f_0,g_0)$.
  The map $g_0$ induces by construction a surjection $\frg(u)^{(0)}\to\frg(r_0)^{(0)}$ and thus the range of the map induced by  $f'\circ f_0\in\Hom_\Arc(\underline r_0, (X,\theta))$ from $\frg(r_0)^{(0)}$ to $G(X,\theta)^{(0)}= X/\theta$ is contained in $Y=\phi(\frg(u)^{(0)}) \subset G(X,\theta)^{(0)}$. Moreover $Y=f(\frg(r)^{(0)})$ and the map $f^{(0)}$ gives an increasing bijection of  $\{0,\ldots r\}$ with $Y$ by construction. Thus, there exists $\rho \in \Hom_\Arc(\underline r_0,\underline r)$ such that $f'\circ f_0=f\circ \rho$ and one gets
 $$
(f',g')\sim(f'\circ f_0,g_0)=(f\circ \rho,g_0)\sim (f, G(\rho)\circ g_0)=(f,g)
$$
which gives the required equivalence.

We have shown that the map $\gamma:p^*(X)\to \Hom_\frg(Z, G(X,\theta))$ is bijective, it follows that $j(p)$ is given by  the flat functor $F:\lbto\longrightarrow \Se$, $F(\underline u)= \Hom_\frg(\frg(u), G(X,\theta))$.  Translating this fact in terms of semifields and  \mods~one obtains the required result.\qed

%Subsection: Sections \iota_n of Mod...
\subsection{Sections $\iota_n$ of $\mmod$ and their action on points}\label{sectgeomrel}

Let $n\geq 0$ and consider, as above,  the archimedean set $\underline n:=(\Z,t^{n+1})$ given by the translation $x\mapsto x+n+1$ on $\Z$. The map 
$x\mapsto kx$ on $\Z$ defines, for any $k\in \nt$, an element $\iota_n(k)\in \Hom_{\Arc\ltimes \nt}(\underline n,\underline n)$ and one obtains a 
homomorphism of semigroups $\iota_n:\nt\to \Hom_{\Arc\ltimes \nt}(\underline n,\underline n)$. Moreover by construction one has $\Mod(\iota_n(k))=k$ $\forall k\in \nt$. This defines, 
for any $n\geq 0$,  a section $\iota_n:\nt \longrightarrow \lbto$ of the functor $\mmod$. More precisely, $\iota_n$ is associated to the functor  which maps   the only object $\bullet$ of $\nt$ to $\underline n$ and sends the morphism $k\in \End_{\nt}(\bullet)$ to $\iota_n(k)\in\Hom_\lbto(\underline n, \underline n)=\End_{\frg}(\frg(n))$.
One has ${\rm Mod}\circ \iota_n=\id$.

\begin{prop}\label{extofscal} Let $p$ be the point of the topos $\widehat{ \N^\times}$ associated to the ordered group $(H,H_+)$. Let $K=H_{\rm max}$ be the semifield associated to $(H,H_+)$. Then the image $\iota_n(p)$ by the section $\iota_n$ is the point of $\epi$ associated by Theorem \ref{main} to the pair $(K,K^{(n+1)})$,
  where for $m\geq 1$, $K^{(m)}$ denotes the \mod~over $K$ obtained from $K$ by extension of scalars using the endomorphism $\fr_m\in \End(K)$, $\fr_m(x):=x^m$ $\forall x\in K$. 
\end{prop}
\proof  By \eqref{yoneda1}, $\iota_n(p)$ is given by  the flat functor $F:\lbto\longrightarrow \Se$
\begin{equation}\label{yoneda3}
F(Z)=p^*(X), \   \   X:\nt\longrightarrow \Se, \  \  X(\bullet):=\Hom_{\lbto}(\iota_n(\bullet),Z).
\end{equation}
One has $\Hom_{\lbto}(\iota_n(\bullet),Z)=\Hom_{\frg}(Z, \frg(n))$. For any $k\in \nt=\End_{\nt}(\bullet)$ the map $X(k):X(\bullet)\to X(\bullet)$ is obtained by composition with $\iota_n(k)\in \End_{\frg}(\frg(n))$.
The inverse image part $p^*$ of the geometric morphism associated to the point $p$ is the tensor product which associates to any $\N^\times$-space $X$ the set
$$
p^*(X):= X\times_{\N^\times} H_+.
$$
For $X=\Hom_{\frg}(Z, \frg(n))$ one gets $p^*(X)=\Hom_{\frg}(Z, G(H^{(n+1)},H))$ where, for $m\geq 1$, $H^{(m)}$ is the totally  ordered set $H$ on which the ordered group $H$ acts by 
$
(h,x)\mapsto x+mh, \forall x\in H^{(m)}, h\in H
$. 
The flat functor $F$ is thus given by $F(Z)=\Hom_{\frg}(Z, G(H^{(n+1)},H))$.
By formulating this construction in terms of semifields and  \mods~one obtains the required result.\qed

%Subsection: Automorphisms of the cyclic and epicyclic categories

\newpage

%Appendix A

\appendix
\section{The action of the barycentric subdivision on  points of $\hat \Delta$}\label{appbary}

We recall that the simplicial category $\Delta$ is the small category with objects the totally ordered sets $[n] :=\{0,\ldots,n\}$,  for each integer $n\geq 0$, and morphisms non-decreasing maps.\newline
In this section we study, using the formalism of topoi, the barycentric subdivision functors $\sd_k: \Delta\longrightarrow\Delta$, for $k\in \nt$ and their action on the points of the simplicial topos $\hat\Delta$. 
\subsection{The barycentric subdivision functors $\sd_k$}
Let $F$ be a  finite, totally ordered set and $k\in \nt$ a positive integer. We define the set 
\begin{equation}\label{sdkf00}
\sd_k(F):=\{0,\ldots,k-1\} \times F  
\end{equation}
to be the cartesian product of the finite ordered set $\{0,\ldots,k-1\}$ with $F$, endowed with the lexicographic ordering. For $f\in \Hom_\Delta(F,F')$ a non-decreasing map (of finite, totally ordered sets), we let 
\begin{equation}\label{sdkf0}
\sd_k(f):= \id \times f: \sd_k(F)\to \sd_k(F')
\end{equation} 
\begin{prop}\label{defnsubdiv} For each $k\in \nt$, \eqref{sdkf00} and \eqref{sdkf0} define an endofunctor $\sd_k: \Delta\longrightarrow\Delta$. They fulfill the property
$$
\sd_{kk'}=\sd_k \circ \sd_{k'}\qqq k,k'\in \nt.
$$
\end{prop}
\proof The totally ordered sets  $\sd_k([n])$ and $[k(n+1)-1]$ have the same cardinality and are canonically isomorphic.
The unique increasing bijection $ \sd_k([n])\to[k(n+1)-1]$ is given by 
$$
(a,i)\mapsto i+a(n+1)\qqq a\in \{0,\ldots,k-1\}, \  i\in \{0,\ldots,n\}.
$$ Let $f\in  \Hom_\Delta([n],[m])$ then by definition $\sd_k(f)\in  \Hom_\Delta(\sd_k([n]),\sd_k([m]))$ is given by
\begin{equation}\label{sdkf}
\sd_k(f)(i+a(n+1))=f(i)+a(m+1)\qqq i,a, \  0\leq i\leq n,~ 0\leq a\leq k-1.
\end{equation}
We check that  $\sd_{kk'}=\sd_k \circ \sd_{k'}$. One has $\sd_{k'}(f)\in  \Hom_\Delta(\sd_{k'}([n]),\sd_{k'}([m]))$ and for any $u\in \{0,\ldots,k'(n+1)-1\}$ and $a\in \{0,\ldots,k-1\}$ one gets 
$$
\sd_k(\sd_{k'}(f))(u+a(k'(n+1)))=\sd_{k'}(f)(u)+a k'(m+1).
$$
Moreover, with $u=i+a'(n+1)$, $0\leq i\leq n, 0\leq a'\leq k'-1$, one has
$$
\sd_{k'}(f)(u)+a k'(m+1)=f(i)+a'(m+1)+a k'(m+1)=f(i)+b(m+1)
$$
where $b=a k'+a'$. Thus for any $b\in \{0,\ldots, kk'-1\}$ one derives
$$
\sd_k(\sd_{k'}(f))(i+b(n+1))=f(i)+b(m+1)
$$
which shows that $\sd_k(\sd_{k'}(f))=\sd_{kk'}(f)$. \endproof

We transfer the functors $\sd_k$ to the opposite category $\dop$ of finite intervals. Recall that by definition, an interval  $I$ is a totally ordered set with a smallest element $b$ and a largest elements  $t\neq b$. The morphisms between intervals are the non-decreasing maps respecting $b$ and $t$, \ie  $f:I\to J$, $f(b_I)=b_J$, $f(t_I)=t_J$. \newline
For all $n\geq 0$ we denote by $n^*:=\{0,\ldots, n+1\}$. The interval $n^*$ parametrizes the hereditary subsets of $[n]$: indeed, to $j\in n^*$ corresponds $[j,n]:=\{x\in [n]\mid x\geq j\}$,  the latter set is empty for $j=n+1$. The duality between $\Delta$ and $\dop$ is then provided by the contravariant functor $\Delta\stackrel{\sim}{\longrightarrow}\dop$,   $[n]\mapsto n^*$, which acts on morphisms as follows
\begin{equation}\label{dualdop}
\Hom_\Delta([n],[m])\ni f  \to f^*\in \Hom_\dop(m^*,n^*)
, \ 
f^{-1}([j,m])=[f^*(j),n] \qqq j\in m^*
\end{equation}
Let  $I$ be an interval, and $k\in \nt$, then one lets $\sdd_k(I)$ to be the quotient of the totally ordered set
$\{0,\ldots,k-1\} \times I= \sd_k(I)$ (with lexicographic ordering) by the equivalence relation
 $(j,t_I)\sim (j+1,b_I)$ for $j\in \{0,\ldots,k-2\}$. This defines an endofunctor $\sdd_k$ of the category of intervals whose action on morphisms sends 
 $f:I\to J$ to $\sdd_k(f)=\id \times f$. By restriction to finite intervals one obtains an endofunctor $\sdd_k: \dop\longrightarrow\dop$. \newline
In particular,  the interval $\sdd_k(n^*)$ has $k(n+2)-(k-1)=k(n+1)+1$ elements and one obtains a canonical  identification of  $\sdd_k(n^*)$  with the hereditary subsets of $\sd_k([n])$ as follows
$$
\{0,\ldots,k-1\} \times n^*\ni (b,j)\mapsto \{ (a,i)\in \{0,\ldots,k-1\} \times [n]\mid a>b \ \text{or}\  a=b\  \&\  j\geq i\}. 
$$
Note that the right hand side of the above formula  depends only upon the class of $(b,j)\in \sdd_k(n^*)$.
\begin{lem} For $f\in \Hom_\Delta([n],[m])$, one has
$(\sd_{k}(f))^*=\sdd_k(f^*).
$
\end{lem}
\proof
The morphism $(\sd_{k}(f))^*$ is defined by the equivalence
$$
\sd_k(f)(x)\geq y\iff (\sd_{k}(f))^*(y)\leq x
$$
Let $x=i+a(n+1)$, $y= j+b(m+1)$ with $0\leq i\leq n$, $0\leq j\leq m$, $0\leq a\leq k-1$, $0\leq b\leq k-1$. Then by \eqref{sdkf} one has $\sd_k(f)(x)=f(i)+a(m+1)$,  thus the condition $\sd_k(f)(x)\geq y$ determines $(\sd_{k}(f))^*(y)$ as follows
$$
f(i)+a(m+1)\geq  j+b(m+1)\iff a>b\ \text{or} \ a=b \  \text {and} \  f(i)\geq j
$$
$$
\iff i +a(n+1)\geq f^*(j)+b(n+1)=\sdd_k(f^*)(y).
$$
This provides the required equality $(\sd_{k}(f))^*=\sdd_k(f^*)
$. \qed

\begin{thm}: The action $\hat \sd_k$ of the geometric morphism $\sd_k$ ($k\in\nt$) on the  points of the topos $\hat \Delta$
is described by the endofunctor $\sdd_k$ on the category of intervals.
\end{thm}
\proof One can prove this theorem using the fact that any point of $\hat \Delta$ is obtained as a filtering colimit of the points associated to the Yoneda embedding of $\dop$ in the category of points of $\hat \Delta$. One shows that on such points the action $\hat \sd_k$ coincides with the functor $\sdd_k:\dop \longrightarrow \dop$. We shall nevertheless find it more instructive to give, in \S \ref{sectexplicit}, a concrete direct proof of the equality between the following two flat functors $\Delta\longrightarrow \Se$ associated to an interval $I$
\begin{equation}\label{twoflats}
F_1([n])=\Hom_\dop(n^*,\sdd_k(I)), \    F_2([n])=\left(\coprod_{m\ge 0} (\Hom(m^*,I)\times  \Hom_\dop(n^*,\sdd_k(m^*)))\right)/\sim
\end{equation}
here  $F_2$ is the  inverse image functor  of the point $p_I$ of $\hat \Delta$ applied to the contravariant functor $Y:\Delta\longrightarrow \Se$, $Y=X\circ\sd_k$, where $X=h_{[n]}$ is the Yoneda embedding, so that
\begin{equation}\label{yy}
Y([m])=X(\sd_k([m]))=\Hom_\Delta(\sd_k([m]),[n])=\Hom_\dop(n^*,\sdd_k(m^*)).
\end{equation}
Thus $F_2$ corresponds to the point $\hat \sd_k(p_I)$ and the equality between $F_1$ and $F_2$ (\cf ~Lemma \ref{lemdescsub})
yields the result. \qed

\begin{cor}: The point of the simplicial topos  $\hat \Delta$ associated to the  interval $[0,1]\subset \R$ is a fixed point for the action of $\N^\times$ on $\hat \Delta$.
\end{cor}
\proof The statement follows using the affine isomorphism
\begin{equation}\label{affineiso}
\sdd_k([0,1])=\{0,\ldots,k-1\}\times [0,1]\to [0,1], \ \ (a,x)\mapsto \frac a k+\frac x k.
\end{equation}

%A.1
\subsection{Canonical decomposition of $\varphi \in \Hom_\dop(n^*,\sdd_k(I))$}\label{sectcandec}
Let $I$ be an interval,  $I/\sim$ be the quotient of $I$  by the identification $b\sim t$. Consider the map  $\pi: \sdd_k(I)\to I/\sim$, $(j,x)\mapsto x$. For
$\varphi \in \Hom_\dop(n^*,\sdd_k(I))$, we define the rank of $\varphi$ as the cardinality of the set $Z=I^o\cap{\rm Range}(\pi \circ \varphi)$, where $I^o=I\setminus \{b,t\}$.
\begin{prop} \label{propapp1} Let 
$\varphi \in \Hom_\dop(n^*,\sdd_k(I))$  and $r$ its rank. Then, one  has a unique decomposition
$$
\varphi=\sdd_k(\alpha)\circ \beta, \ \ \beta\in \Hom_\dop(n^*,\sdd_k(r^*)), \ \alpha\in \Hom_\dop(r^*,I).
$$
Moreover, the morphism $\alpha\in \Hom_\dop(r^*,I)$ is the unique increasing injection $\{1, \ldots, r\}\hookrightarrow I^o$ which admits  $Z$ as range. The composite $\pi_r\circ \beta$ is surjective, where $\pi_r: \sdd_k(r^*)\to r^*/\sim$ is the canonical surjection.
\end{prop}
\proof Let $\alpha\in \Hom_\dop(r^*,I)$ be the map  whose restriction to $\{1, \ldots, r\}$ is the unique increasing injection  into $I^o$ which admits  $Z$ as range. Recall that an element $x\in \sdd_k(I)$ is given by a pair 
$x=(j,y)\in\{0,\ldots,k-1\} \times I= \sd_k(I)$ with the identifications $(j,t)\sim (j+1,b)$ for $j\in \{0,\ldots,k-2\}$. Similarly an element $z\in \sdd_k(r^*)$ is given by a pair $z=(i,u)\in\{0,\ldots,k-1\} \times r^*= \sd_k(r^*)$ with the identifications $(j,r+1)\sim (j+1,0)$ for $j\in \{0,\ldots,k-2\}$. Let $s\in n^*=\{0,\ldots,n+1\}$, then $\varphi(s)\in \sdd_k(I)$ is given by a pair 
$\varphi(s)=(j,y)\in\{0,\ldots,k-1\} \times I= \sd_k(I)$ unique up to the above identifications. If $y\in \{b,t\}$, one defines $\beta(s)=(j,0)\in \sdd_k(r^*)$ if $y=b$, and $\beta(s):=(j,r+1)\in \sdd_k(r^*)$ if $y=t$. This definition is compatible with the identifications. Let us now assume that $y\notin \{b,t\}$. Then $y\in Z$ and there exists a unique element $v\in\{1, \ldots, r\}$ such that $y=\alpha(v)$. One then defines $\beta(s):=(j,v)\in\sdd_k(r^*)$. The map $\beta : n^*\to\sdd_k(r^*)$ so defined is non-decreasing, \ie for $s<s'$ one has $\beta(s')\geq \beta(s)$ since the inequality $\varphi(s')\geq \varphi(s)$ shows that either $j'>j$ in which case $(j',v')\geq (j,v)$ is automatic, or $j=j'$ and in that case $y'>y$ which shows that $v'\geq v$. Moreover since $\sdd_k(\alpha)=\id \times \alpha$ one has $\varphi=\sdd_k(\alpha)\circ \beta$. \newline
We prove the uniqueness of this decomposition. Since $\varphi \in \Hom_\dop(n^*,\sdd_k(I))$ preserves the base points, Range $(\pi \circ \varphi)$ contains the base point and its cardinality is $r+1$. Thus the  map $\alpha\in \Hom_\dop(r^*,I)$ is the unique map  whose restriction to $\{1, \ldots, r\}$ is the  increasing injection  to $I^o$ and which admits  $Z=I^o\cap$ Range $(\pi \circ \varphi)$ as range. Moreover $\alpha$ is injective and so is $\sdd_k(\alpha)$. Thus the map $\beta\in \Hom_\dop(n^*,\sdd_k(r^*))$ is uniquely determined by the equality $\varphi=\sdd_k(\alpha)\circ \beta$. Finally  $\pi_r\circ \beta:n^*\to r^*/\sim$ is surjective since otherwise the range of $\sdd_k(\alpha)\circ \beta$ would be strictly smaller than the range of $\varphi$. \qed

\begin{cor}\label{surj}
For any interval $I$ the map 
$$
\Hom_\dop(n^*,I)\times \Hom_\dop(n^*,\sdd_k(n^*))\to \Hom_\dop(n^*,\sdd_k(I)),\quad (\alpha,\beta)\mapsto \sdd_k(\alpha)\circ \beta
$$
is surjective.
\end{cor}

%A2
\subsection{Explicit description of the isomorphism $F_1\simeq F_2$}\label{sectexplicit}

Let $F_j: \Delta\longrightarrow \Se$ be the flat functors  defined in \eqref{twoflats}. By definition 
\begin{equation}\label{explicit}
  F_2([n])=\left(\coprod_{m\ge 0} (\Hom(m^*,I)\times  \Hom_\dop(n^*,\sdd_k(m^*)))\right)/\sim
\end{equation}
where the equivalence relation is generated by  
$$
(\alpha\circ f,\beta)\sim (\alpha,\sdd_k(f)\circ \beta)
$$
for $f\in \Hom_\dop(m^*,r^*)$, $\beta\in  \Hom_\dop(n^*,\sdd_k(m^*))$, $\alpha\in \Hom(r^*,I)$.

\begin{lem} \label{lemdescsub}The map $\Phi: F_2([n])\to F_1([n])=\Hom_\dop(n^*,\sdd_k(I))$, $(\alpha,\beta)\mapsto \sdd_k(\alpha)\circ \beta$ is a  bijection of sets.
\end{lem}
\proof
The map $\Phi$ is well defined since $\Phi(\alpha\circ f,\beta)=\Phi(\alpha,\sdd_k(f)\circ \beta)$.
Corollary \ref{surj} shows that $\Phi$ is surjective. To show the injectivity it is enough to prove that for any
$(\alpha,\beta)\in \Hom(m^*,I)\times  \Hom_\dop(n^*,\sdd_k(m^*))$ one has $(\alpha,\beta)\sim (\alpha_c,\beta_c)$ where
$$
\varphi=\sdd_k(\alpha_c)\circ \beta_c, \ \beta_c \in \Hom_\dop(n^*,\sdd_k(r^*)), \ \alpha_c\in \Hom_\dop(r^*,I)
$$
is the canonical decomposition of $\varphi=\sdd_k(\alpha)\circ \beta$. 

One has the canonical decomposition $\beta=\sdd_k(\alpha_0)\circ \beta_0$ with 
$\beta_0\in \Hom_\dop(n^*,\sdd_k(\ell^*))$, $\ps_\ell\circ \beta_0$ surjective. Thus $(\alpha,\beta)\sim (\alpha\circ \alpha_0,\beta_0)$. Since $\ps_\ell\circ \beta_0$ is surjective, Range $\alpha\circ \alpha_0 \subset $ Range $\ps \circ \varphi$ = Range $\alpha_c$, and thus  $\alpha\circ \alpha_0=\alpha_c\circ \rho$
$$
(\alpha,\beta)\sim(\alpha\circ \alpha_0,\beta_0)=(\alpha_c\circ \rho,\beta_0)\sim (\alpha_c, \sdd_k(\rho)\circ \beta_0)=(\alpha_c,\beta_c).
$$\endproof

%A3

\subsection{The small category $\dop\ltimes \N^\times$}\label{sectcrossprod}

We denote by $\dop\ltimes \N^\times$  the small category  semi-direct product of $\dop$ by the action of $\nt$ implemented by the endofunctors $\sdd_k$, for $k\in \nt$. It has the same objects as  $\dop$ while one adjoins to the collection of morphisms of $\dop$ the new morphisms  $\pi_n^k:\sdd_k(n^*)=\left(k(n+1)-1\right)^*\to n^*$ such that
\begin{equation}\label{complawpink}
\pi_n^k\circ \pi_{k(n+1)-1}^\ell=\pi_n^{k\ell}\in \Hom_{\dop\ltimes \N^\times}(\left(k\ell(n+1)-1\right)^*,n^*)
\end{equation}
where  $\pi_n^k$ implements the endofunctor $\sdd_k$, \ie 
\begin{equation}\label{complawcross}
\alpha\circ \pi_n^k= \pi_m^k\circ \sdd_k(\alpha)\qqq \alpha \in \Hom_\dop(n^*,m^*).
\end{equation}
Using this set-up one checks that any morphism $\phi$ in $\dop\ltimes \N^\times$ is uniquely of the form 
$\phi=\pi_n^k\circ \alpha$ with $\alpha$ a morphism in $\dop$. Any such $\phi$ compose as follows
\begin{equation}\label{complacomp}
(\pi_m^k\circ\beta)\circ (\pi_n^\ell\circ \alpha)=\pi_m^{k\ell}\circ(\sdd_\ell(\beta)\circ \alpha)
\end{equation}
where $\alpha\in \Hom_\dop(r^*,\left(\ell(n+1)-1\right)^*)$, $\sdd_\ell(\beta)\in\Hom_\dop(\left(\ell(n+1)-1\right)^*,\left(k\ell(m+1)-1\right)^*)$ so that $\sdd_\ell(\beta)\circ \alpha$ makes sense and belongs to $\Hom_\dop(r^*,\left(k\ell(m+1)-1\right)^*)$.
Using Proposition \ref{defnsubdiv} one checks that, if one takes \eqref{complacomp} as a definition, the product is associative.\newline

Note that the association $\phi=\pi_n^k\circ \alpha\mapsto k\in \N^\times$ determines a functor $\dop\ltimes \N^\times\longrightarrow \N^\times$ which is not injective on objects (unlike the inclusion $\dop\longrightarrow \dop\ltimes \N^\times$).

Let $\fin_*$ be the category of finite pointed sets and let $\cF$ be the functor which associates to an interval $I$ the pointed set $I_*=I/\sim$ with base point the  class of $b\sim t$. To any morphism of intervals $f:I\to J$ corresponds the quotient map $f_*$ which preserves the base point. By restricting $\cF$ to $\dop$ one gets a covariant functor $\cF:\dop\longrightarrow \fin_*$. The following Proposition shows that $\cF$ can be extended to $\dop\ltimes \N^\times$.

\begin{prop}: For any $n\geq 0,k\in\N^\times$, let $(\pi_n^k)_*:\cF(\sdd_k(n^*))\to\cF(n^*)$ be given by the residue modulo $n+1$. Then the  extension of the functor  $\cF$  on morphisms given by
$$
\phi=\pi_n^k\circ \alpha\mapsto \phi_*:=(\pi_n^k)_*\circ \alpha_*
$$
 determines a functor $\cF:\dop\ltimes \N^\times\longrightarrow \fin_*$.
\end{prop}
\proof One checks directly that the definition of $(\pi_n^k)_*$ is compatible with the rules \eqref{complawpink} and \eqref{complawcross} so that the required functoriality follows. \qed

%A4
\subsection{The subset $\Sigma_n^k\subset \Hom_\dop(n^*,\sdd_k(n^*))$}

For each $n\geq 0$ and $k\in \nt$, we introduce the following set 
\begin{equation}\label{defnsigmank}
\Sigma_n^k:=\{\alpha \in \Hom_\dop(n^*,\sdd_k(n^*))\mid (\pi_n^k)_*\circ \alpha_*~\text{is surjective}\}\subset \Hom_\dop(n^*,\sdd_k(n^*)).
\end{equation}
Here $(\pi_n^k)_*\circ \alpha_*: \{0, \ldots ,n\}\to\{0, \ldots ,n\}$ preserves the base point $0$, moreover it is surjective if and only if it is a permutation of $\{0, \ldots ,n\}$ fixing $0$. One thus obtains a canonical map  to the group $S_n$ of permutations of $\{1, \ldots ,n\}$ defined as follows
\begin{equation}\label{maptoperm}
\pee:\Sigma_n^k\to S_n, \  \  \pee(\alpha):=(\pi_n^k)_*\circ \alpha_*.
\end{equation}
One gets a parametrization of the set $\Sigma_n^k$ as follows, let $f:\{1, \ldots ,n\}\to \{0, \ldots ,k-1\}$ be an arbitrary map of sets  and consider the following subset with $n$ elements 
 \begin{equation}\label{setXf}
X_f:=\{j+(n+1)f(j)\mid j\in \{1, \ldots ,n\}\}\subset \{1, \ldots ,k(n+1)-1\}
\end{equation}
Then, as shown by the next Lemma \ref{lemparam}, there exists a unique element $s(f)\in \Hom_\dop(n^*,\sdd_k(n^*))$ whose range is $X_f$. It is obtained by labeling the elements of $X_f$ in the lexicographic order, with $s(f)(0)=0$ and $s(f)(n+1)=k(n+1)$. 
\begin{lem}\label{lemparam} The map $f\mapsto s(f)$ is a bijection of the set of all maps $\{1, \ldots ,n\}\to \{0, \ldots ,k-1\}$ with $\Sigma_n^k$.
\end{lem}
\proof First the range of $s(f)$ when taken modulo $n+1$ contains all $ j\in \{1, \ldots ,n\}$ and thus $s(f) \in \Sigma_n^k$. By construction the map $f\mapsto s(f)$ is injective since $f(j)$ is the only element of $X_f$ which is congruent to $j$ modulo $n+1$. Let us show that it is surjective. Let $\alpha \in \Sigma_n^k$, then since $(\pi_n^k)_*\circ \alpha_*$ is surjective there exists for each
$ j\in \{1, \ldots ,n\}$ an $f(j)\in \{0, \ldots ,k-1\}$ such that $j+f(j)(n+1)\in {\rm Range}(\alpha)$. One has ${\rm Range}(\alpha)=X_f$ and it follows that $\alpha=s(f)$.\qed

\begin{prop}: \label{propdesc}The range of  $\pee:\Sigma_n^k\to S_n$ is the set of permutations whose descent number is $d=k-1-v< k$. Moreover
$$
\#\{\alpha \in \Sigma_n^k\mid \pee(\alpha)=\sigma\}={n+v\choose n}.
$$
\end{prop}
\proof Given a permutation $\sigma\in S_n$, a lift of $\sigma$ in $\Sigma_n^k$ is of the form
 \begin{equation}\label{desc0}
\alpha(j)=\sigma(j)+k(j)(n+1)
\end{equation}
where $k:\{1, \ldots ,n\}\to \{0, \ldots ,k-1\}$ is non-decreasing and such that 
 \begin{equation}\label{desc}
\sigma(j+1)<\sigma(j)\implies k(j+1)>k(j)
\end{equation}
By definition (\cf\cite{Loday}) the descent number number $d$ of $\sigma$ is the number of $ j\in \{1, \ldots ,n-1\}$ such that $\sigma(j+1)<\sigma(j)$. Thus \eqref{desc} implies that $k(n)-k(1)\geq d$ and one gets $k-1\geq d$. Conversely assume that the descent number number $d$ of $\sigma$ is $d=k-1-v$ with $v\geq 0$. Then the lifts $\alpha$ of $\sigma$ correspond to sequences $\epsilon(j)\geq 0$, $ j\in \{1, \ldots ,n\}$,  such that 
$$
\sum \epsilon(j)\leq k-1, \ \  \sigma(j+1)<\sigma(j)\implies \epsilon(j)>0.
$$
Given such a sequence one lets $k(j)=\sum_1^j\epsilon(i)$ and one derives a lift of $\sigma$ using \eqref{desc0}. Replacing $\epsilon(j)$ by $\epsilon(j)-1$ when $\sigma(j+1)<\sigma(j)$ one gets that the number of lifts is the number of sequences $\epsilon'(j)\geq 0$,  $ j\in \{0, \ldots ,n\}$ such that $\sum \epsilon'(j)=v$, \ie $n+v\choose n$. \qed

%A5

\subsection{The $\lambda$-operations as elements of  $\Z[\dop\ltimes \N^\times]$}\label{A6}
Let  $\Z[\dop]$ be the ring of finite formal sums of elements of $\dop$ with integral coefficients,
where the product law is given by the composition of morphisms in $\dop$ whenever they are composable, otherwise is defined to be  $0$.
One first introduces the elements
 \begin{equation}\label{sigmank}
\sigma_n^k:=\sum_{\alpha\in\Sigma_n^k}\epsilon(\pee(\alpha)) \alpha\in \Z[\dop].
\end{equation}
Then, by composing with the maps $\pi_n^k$ one obtains the following elements of 
the ring $\Z[\dop\ltimes \N^\times]$
 \begin{equation}\label{lambdank}
\Lambda_n^k:=\sum_{\alpha\in\Sigma_n^k}\epsilon(\pee(\alpha)) \pi_n^k\circ \alpha.
\end{equation}
Next, we show that these elements $\Lambda_n^k\in \Z[\dop\ltimes \N^\times]$ fulfill the law
 \begin{equation}\label{lambdank1}
\Lambda_n^k \Lambda_n^\ell=\Lambda_n^{k\ell} \qqq k, \ell, n\in \N.
\end{equation}
One has $\Lambda_n^k \Lambda_n^\ell=\pi_n^k \sigma_n^k \pi_n^\ell \sigma_n^\ell=\pi_n^{k\ell}\sdd_\ell(\sigma_n^k)\sigma_n^\ell$
where $\sdd_\ell$ is extended to $\Z[\dop]$ by linearity. Thus \eqref{lambdank1} follows from the next Lemma \ref{lambdank2} which shows that the  map 
$$
\Sigma_n^k\times \Sigma_n^\ell\to \Sigma_n^{k\ell}, \  \  (\alpha,\beta)\mapsto \sdd_\ell(\alpha)\beta
$$
is bijective 
and that the signatures are multiplicative.

\begin{lem}\label{lambdank2}
$(i)$~Let $\alpha \in \Sigma_n^k$ and $\beta\in  \Sigma_n^\ell$. Then $\sdd_\ell(\alpha)\beta\in \Sigma_n^{k\ell}$ and~
  $\pee(\sdd_\ell(\alpha)\beta)=\pee(\alpha)\circ\pee(\beta)$.
$(ii)$~Let $f:\{1, \ldots , n\}\to \{0, \ldots , k-1\}$,
$g:\{1, \ldots , n\}\to \{0, \ldots , \ell-1\}$ be arbitrary maps of sets. Then one has
$$
\sdd_\ell(s(f))\circ s(g)=s(h)
$$
where $h:\{1, \ldots , n\}\to \{0, \ldots , k\ell -1\}$ is given by
$$
h(j)=f(j)+kg(\sigma^{-1}(j)), \ \  \sigma=\pee(s(f)).
$$
\end{lem}
\proof $(i)$~One has $\alpha \in \Hom_\dop(n^*,\sdd_k(n^*))$, $\sdd_\ell(\alpha)\in\Hom_\dop(\sdd_\ell(n^*),\sdd_{\ell k}(n^*))$ and
$\beta \in \Hom_\dop(n^*,\sdd_\ell(n^*))$. Thus $\sdd_\ell(\alpha)\beta\in \Hom_\dop(n^*,\sdd_{\ell k}(n^*))$. Moreover, by applying the functor $\cF:\dop\ltimes N^\times\to \fin_*$ one derives 
$$
\cF(\pi_n^{k\ell}(\sdd_\ell(\alpha)\beta)=\cF(\pi_n^k\alpha\circ\pi_n^\ell\beta)=\cF(\pi_n^k\alpha)\cF(\pi_n^\ell\beta)=\pee(\alpha)\circ\pee(\beta).
$$
This shows that $\sdd_\ell(\alpha)\beta\in \Sigma_n^{k\ell}$ and
  $\pee(\sdd_\ell(\alpha)\beta)=\pee(\alpha)\circ\pee(\beta)$.

$(ii)$~We determine the range of $\gamma=\sdd_\ell(s(f))\circ s(g)$. It is the image by $\sdd_\ell(s(f))$ of the range of $s(g)$ which is by construction $X_g=\{j+g(j)(n+1)\mid j\in \{1, \ldots , n\}\}$ (ignoring the base points). Now one has 
$$
\sdd_\ell(s(f))(i+a(n+1))=s(f)(i)+ka(n+1), 
$$
and thus
$$
\sdd_\ell(s(f))(j+g(j)(n+1))=s(f)(j)+kg(j)(n+1).
$$
Let $\sigma=\pee(s(f))$, then one has $s(f)(j)=\sigma(j)+f(\sigma(j))(n+1)$ and thus taking $j=\sigma^{-1}(i)$ one gets
$s(f)(\sigma^{-1}(i))=i+f(i)(n+1)$ and 
$$
\sdd_\ell(s(f))(j+g(j)(n+1))=i+f(i)(n+1)+kg(\sigma^{-1}(i))(n+1)
$$
which determines uniquely the range of $\gamma=\sdd_\ell(s(f))\circ s(g)$ in the required form. \qed

%A6
\subsection{The commutation $\Lambda_{n-1}^k b=b \Lambda_n^k$}\label{A7}

We define the Hochschild boundary operator $b$,  with one component $b_n$ for each $n\geq 1$, as the element 
 \begin{equation}\label{bn}
b_n:=\sum_0^n(-1)^id_i^n\in \Z[\dop],
\end{equation}
where for each $n\geq 1$ and $i \in \{0,\ldots,n\}$ one lets
$$
d_i^n\in \Hom_\dop(n^*,(n-1)^*), \   d_i^n(j):=\left\{
                 \begin{array}{ll}
                  j, & \hbox{if $j\leq i$;} \\
                   j-1, & \hbox{if $j>i$.} \\
                                   \end{array}
               \right.
$$
One has 
$
d_i\circ d_j= d_{j-1}\circ d_i\qqq i<j
$, and this shows that $b_n\circ b_{n+1}=0$. Next, we provide the geometric proof (\!\!\cite{MCarthy}) of the commutation of the $\lambda$-operations with the Hochschild boundary.
\begin{prop}\label{propblamcom} The elements $\Lambda_n^k\in \Z[\dop\ltimes \N^\times]$ commute with the Hochschild boundary
$
b
$
\ie one has $\Lambda_{n-1}^k b_n=b_n \Lambda_n^k$, $\forall n\geq 1$.
\end{prop}
\proof Note first that since the number of terms in the formula \eqref{sigmank} defining $\Lambda_n^k$ is the cardinality $k^n$ of $\Sigma_n^k$, the expansion of the sum 
$$
\Lambda_{n-1}^k b_n=\sum_{\delta \in \Sigma_{n-1}^k\atop 0\leq j\leq n} (-1)^j \epsilon(\pee(\delta))\, \pi_{n-1}^k\circ \delta\circ d_j^n
$$ 
contains $(n+1)k^{n-1}$ terms while the expansion of 
$$
b_n \Lambda_n^k=\sum_{\alpha \in \Sigma_{n}^k\atop 0\leq j\leq n} (-1)^i \epsilon(\pee(\alpha))\, d_i^n\circ\pi_{n}^k\circ \alpha
$$ 
contains $(n+1)k^n$ terms.  One has 
$$
d_i^n\circ\pi_{n}^k=\pi_{n-1}^k\circ\sdd_k(d_i^n)
$$
and thus the equality to be proved is 
$$
\sum_{\delta \in \Sigma_{n-1}^k\atop 0\leq j\leq n} (-1)^j \epsilon(\pee(\delta))\,  \delta\circ d_j^n=
\sum_{\alpha \in \Sigma_{n}^k\atop 0\leq i\leq n} (-1)^i \epsilon(\pee(\alpha))\,\sdd_k(d_i^n)\circ \alpha.
$$

We show geometrically that among the $(n+1)k^n$ terms on the right hand side of the above formula there are $(n+1)k^{n-1}$ terms which correspond to the left hand side while the others cancel in pairs. Indeed, we construct a natural correspondence (\cf~ \eqref{triang2} below) between the terms $\sdd_k(d_i^n)\circ \alpha$ and the faces of the simplices in the triangulation of the standard simplex $\underline\Delta^n=\{(x_b)\mid 0\leq x_1\leq \cdots\leq x_n\leq 1\}$ given by the barycentric subdivision into $k^n$ simplices $\underline\Delta(\alpha)$ parametrized by $\alpha \in \Sigma_{n}^k$. The faces which belong to the interior of $\underline\Delta^n$ appear twice and with opposite orientations. This corresponds to the cancelation by pairs of $(n+1)(k^n-k^{n-1})$ terms. On the other hand, the faces which belong to the boundary of $\underline\Delta^n$ correspond to the remaining $(n+1)k^{n-1}$ terms which correspond to the terms on the left hand side.
\begin{figure}
\begin{center}
\includegraphics[scale=0.5]{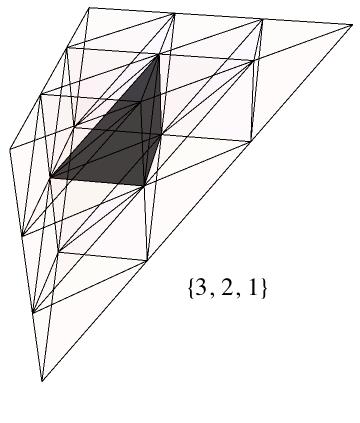}
\end{center}
\caption{The cell of the barycentric subdivision for $n=3, k=3$ corresponding to the permutation with descent number $2$.\label{clock3} }
\end{figure}

Consider the affine isomorphism \eqref{affineiso}, and  associate to  the elements $x\in \Hom_\dop(n^*,[0,1])=\underline\Delta^n$   the map of intervals 
$$
\sdd_k(x)\in \Hom_\dop(\sdd_k(n^*),\sdd_k([0,1]))\cong \Hom_\dop(\sdd_k(n^*),[0,1]).
$$
In terms of the coordinates $x(j)$ one has
 \begin{equation}\label{triang0}
\sdd_k(x)(j)=\frac a k+\frac{ x(b)} {k} \qqq j=b+a(n+1)\in \sdd_k(n^*).
\end{equation}
One obtains in this way  a triangulation of $\underline\Delta^n$ by the $k^n$ simplices 
 \begin{equation}\label{triang}
\underline\Delta(\alpha)=\{\sdd_k(x)\circ \alpha\mid x\in \underline\Delta^n\}\qqq \alpha \in \Sigma_n^k.
\end{equation}
To each simplex $\underline\Delta(\alpha)$ corresponds a permutation $\pee(\alpha)=(\pi_n^k\circ\alpha)_*$. Geometrically  the map $x\mapsto \sdd_k(x)\circ \alpha$ is of the form 
$
(x_b)\mapsto (y_b), \  y_b=\frac 1k x_{\sigma(b)}+s_b, \  \sigma=\pee(\alpha)
$.
For each element $\beta\in \Hom_\dop(n^*,\sdd_k((n-1)^*))$ we consider the subset of $\underline\Delta^n$ given by
 \begin{equation}\label{triang2}
F(\beta)=\{\sdd_k(y)\circ \beta\mid y\in \underline\Delta^{n-1}\}, \  \  \underline\Delta^{n-1}=\Hom_\dop((n-1)^*,[0,1])
\end{equation}

\begin{lem}\label{lembb}
$(i)$~Let $\beta,\beta'\in  \Hom_\dop(n^*,\sdd_k((n-1)^*))$ be such that the faces $F(\beta)$, $F(\beta')$ are non-degenerate and equal then $\beta=\beta'$.

$(ii)$~Let $\alpha \in \Sigma_n^k$ and $\underline\Delta(\alpha)$ the corresponding simplex as in \eqref{triang}. Then the faces of 
  $\underline\Delta(\alpha)$ are the $F(\sdd_k(d_i^n)\alpha)$ for $i\in \{0, \ldots , n\}$.
\end{lem}
\proof $(i)$~This follows since the barycenter of a non-degenerate  face $F(\beta)$ determines $\beta$. More precisely since one assumes that  $F(\beta)$ is non-degenerate the map 
$
 \underline\Delta^{n-1}
\to F(\beta)$, $y\mapsto \sdd_k(y)\circ \beta
$
is affine and injective and thus the barycenter of $F(\beta)$ is the image $\sdd_k(z)\circ \beta$ of the barycenter $z\in \Hom_\dop((n-1)^*,[0,1])$ of $ \underline\Delta^{n-1}$. One 
has $z(j)=j/n\in [0,1]$, $\forall j\in \{0, \ldots , n\}$. It follows from \eqref{triang0} that $\sdd_k(z)(i)=i/(nk)$, $\forall i\in \{0, \ldots , nk\}$. Thus the coordinates $\sdd_k(z)\circ \beta(i)$ of the barycenter of $F(\beta)$ determine $\beta$.

$(ii)$~Using \eqref{triang2}, one has
$$
F(\sdd_k(d_i^n)\alpha)=\{\sdd_k(y)\circ \sdd_k(d_i^n)\alpha\mid y\in \underline\Delta^{n-1}\}=\{\sdd_k(y\circ d_i^n)\alpha\mid y\in \underline\Delta^{n-1}\}
$$
and the result follows since the faces of  the standard simplex $\underline\Delta^n=\Hom_\dop(n^*,[0,1])$ are the 
$$
\partial_i \underline\Delta^n=\{y\circ d_i^n\mid y\in \underline\Delta^{n-1}\}\subset \Hom_\dop(n^*,[0,1]).
$$
\qed
%Appendix B
\section{Epicyclic modules and the $\lambda$-operations}\label{appcyclic}

In this appendix we first review the definition of cyclic homology for cyclic modules, then we give a detailed description of the action of the $\lambda$-operations on the $(b,B)$ bicomplex in  the context of epicyclic modules. Since the cyclic category is isomorphic to its dual, we  use a covariant definition of cyclic modules which turns out to be more convenient when one discusses the epicyclic construction. 

We recall that the cyclic category $\Lambda$ is the full subcategory of the category $\Arc$  with one object $\underline n$ for each non-negative integer $n\ge 0$, given by the archimedean set $\underline n=(\Z,\theta)$, $\theta(x)=x+n+1$ $\forall x\in \Z$.  Thus the morphisms between objects are given by $\Hom_\Lambda(\underline n,\underline m) := \cC(n+1,m+1)$, where for each pairs of integers $a,b>0$, $\cC(a,b)$ denotes the set of equivalence classes of maps $f: \Z\to\Z$, with $f(x)\ge f(y)$, $\forall x\ge y$ and $f(x+a) = f(x) + b,~\forall x\in\Z$ (\cf\cite{topos} \S 2.2). The equivalence relation is defined as follows
\[
f\sim g~\Longleftrightarrow~\exists k\in\Z,~g(x) = f(x)+kb\quad\forall x\in\Z.
\]
The cyclic category has a decomposition $\Lambda = \Delta\mathbf C$ which describes $\Lambda$ as an extension of the small simplicial category $\Delta$ by means of a new generator $\tau_n\in \mathbf C_{n+1}:= \Aut_\Lambda(\underline n)$ for each $n\ge 0$, and fulfilling the relations (in terms of the faces $\delta_j$ and degeneracies $\sigma_j$, $2 \le j\le n$, which describe a presentation of $\Delta$: \cf~\opcit for details)
\begin{align*}%\label{pres1}
&&\tau_n^{n+1}=id,&\notag\\
\tau_n\circ \sigma_0&=\sigma_n\circ \tau_{n+1}^2 &&&   \tau_n\circ \sigma_j &=\sigma_{j-1}\circ \tau_{n+1} \qqq j\in \{1, \ldots, n\}\\%\label{pres2}
\tau_n\circ \delta_0 &=\delta_n  &&&  \tau_n\circ \delta_j&=\delta_{j-1}\circ \tau_{n-1} \qqq j\in \{1, \ldots, n\}.
\end{align*}
The epicyclic category is the full subcategory of the category $\arcnt$ of archimedean sets whose objects are the archimedean sets $\underline n$. We recall that for $(X,\theta)$ an archimedean set and $k>0$  an integer, the pair $(X,\theta^k)$ is also an archimedean set which we denote as 
\begin{equation}\label{psik}
\Psi_k(X,\theta):=(X,\theta^k).
\end{equation}

The epicyclic category $\tilde \Lambda$ is obtained (\cf \cite{bu1}, Definition 1.1) by adjoining to $\Lambda$ new morphisms $\ps_n^k:\Psi_k(\underline n)\to \underline n$ for $n\geq 0$, $k\geq 1$, which fulfill the following relations 
\begin{enumerate}
\item $\ps_n^1={\rm id_n} $, $\ps_n^\ell \circ \ps_{\ell(n+1)-1}^{k}=\ps_n^{k\ell}$
\item $\alpha\, \ps_m^k=\ps_n^k\, {\rm Sd_k}(\alpha)$, for any $\alpha\in \Hom_\Delta([m],[n])$
\item $\tau_n \,\ps_n^k=\ps_n^k\, \tau_{k(n+1)-1}$
\end{enumerate}
where $\sd_k: \Delta\longrightarrow\Delta$ is the barycentric subdivision functor (\eqref{sdkf00}, \eqref{sdkf0}).   

The canonical inclusion $\dop\subset \Lambda$ extends to an inclusion $\dop\ltimes \nt\subset \lbt$ which associates to the morphisms $\pi_n^k\in \Hom_{\dop\ltimes \nt}(\sdd_k(n^*), n^*)$ the morphism $\ps_n^k\in \Hom_\lbt(\Psi_k(\underline n), \underline n)$.%B1

\subsection{Cyclic homology and cyclic modules}

\begin{defn}\label{cycmod} A {\em cyclic} module $E$ is a  {\em covariant} functor $\Lambda\longrightarrow\Ab$ from the cyclic category  to the category of abelian groups.
\end{defn}
We briefly recall the construction of the normalized  $(b,B)$-bicomplex of a cyclic module $E$. We keep the notations of \ref{A6} and \ref{A7}. For each integer  $n\geq 0$, one defines  the element $B(n)=B_0A\in \Z[\Lambda]$, for 
$
A:=\sum\epsilon(\pee(\tau))\tau$, where the summation is taken over $ \tau \in \Aut_\Lambda(\underline n)$,
while the map $B_0\in \Hom_\Lambda(\underline n, \underline {n+1})$ is defined as follows
\begin{equation}\label{defnBo} 
 B_0(i+a(n+1))=1+i+a(n+2)\qqq i, 0\leq i\leq n, ~a\in\Z.
\end{equation}
Let  $s_j=s_j^{n}\in \Hom_\dop(n^*,(n+1)^*)$ be the degeneracies, for $0\leq j\leq n$. They are given by 
$$
s_j^{n}(i)=i \, \  \text{if} \, \  i\leq j, \  \  s_j^{n}(i)=i+1 \,\ \text{if}\  \, i> j
$$
so that the image of $s_j^{n}$ does not contain $j+1$. Note that $B_0\neq s_j^{n}$, $\forall j,n$. 
For $n\geq 0$ let
 \begin{equation}\label{chains}
C_n(E):=E(n)/V(n), \ \ V(n):=\oplus \,{\rm Im}(s_j)    
\end{equation}
where the $s_j=s_j^{n-1}\in \Hom_\dop((n-1)^*,n^*)$ are the degeneracies for $0\leq j\leq n-1$. 
Both the operators $b$ and $B$ pass to the quotient and one has (\cf \cite{CoExt}, \cite{Loday})
\begin{lem}: $(i)$~ 
$
B\left(\oplus \,{\rm Im}(s_j)    \right)\subset \left(\oplus \,{\rm Im}(s_\ell)    \right)~\forall j,\ell.
$

$(ii)$~The subspace $W(n):={\rm Im}(B_0)\oplus V(n)\subset E(n)$ is invariant under cyclic permutations: $E(t)W(n)\subset W(n)$.

$(iii)$~$
{\rm Im}(B\circ B_0)\subset \left(\oplus \,{\rm Im}(s_\ell)    \right),
~\forall \ell$.

$(iv)$~ 
$
b\left(\oplus \,{\rm Im}(s_j)    \right)\subset \left(\oplus \,{\rm Im}(s_\ell)    \right),~\forall j,\ell
$.

$(v)$~$b^2=B^2=bB+Bb=0$ on  $C_n(E)$.
\end{lem}

We set the   $(b,B)$-bicomplex in negative degrees as follows
\begin{equation}\label{chains2}
C^{\alpha,\beta}:=C_{\alpha-\beta}(E)\, \  \text{if} \, \ \alpha\geq \beta, \ \ \alpha\leq 0, \  \ C^{\alpha,\beta}:=\{0\}, \  \text{otherwise},
\end{equation}
and one defines the cyclic homology of the cyclic module $E$ by
$$
HC_n(E):=\H^{-n}\left( \bigoplus_{\alpha\leq 0} C^{\alpha,\beta}, b+B\right)
$$

%B2
\subsection{Epicyclic modules}

We generalize the above set-up to the epicyclic category.
\begin{defn} An {\em epicyclic} module $E$ is a  {\em covariant} functor $\tilde\Lambda\longrightarrow\Ab$ from the epicyclic category  to the category of abelian groups.
\end{defn}

The inclusion of categories $\Lambda \longrightarrow \tilde \Lambda$ turns the operator $B$ into an element of the ring $\Z[\lbt]$ for each  $n\ge 0$ (we denote it by $B(n)$). Moreover, the inclusion of categories $\dop\ltimes \N^\times\longrightarrow \tilde \Lambda$ allows one to view the operators $\Lambda_n^k\in\Z[\dop\ltimes \N^\times]$ as elements of the ring $\Z[\lbt]$. In this section we give a proof of the commutation relation $\Lambda_{n+1}^k B =kB  \Lambda_n^k$ (\cf \cite{Loday}, \cite{MCarthy}) at the categorical level.
We use the following two  variants of the sets $\Sigma_n^k$ as in \eqref{defnsigmank}, obtained by
replacing the category $\dop$ with $\Delta$ and $\Lambda$ respectively.  We define, using the canonical inclusion $\Delta\subset\Lambda$
\begin{equation}\label{dnk}
\Delta_n^k=\{\alpha \in \Hom_\Delta([n],\sd_k([n]))\mid (\ps_n^k\circ\alpha)_* \, \text{is bijective}\}
\end{equation}
and 
\begin{equation}\label{cnk}
\Gamma_n^k=\{\alpha \in \Hom_\Lambda(\underline{n},\Psi_k(\underline{n}))\mid (\id_n^k\circ\alpha)_* \, \text{is bijective}\}
\end{equation}
where it follows from \eqref{psik} that $\Psi_k(\underline{n})=(\Z,t^{k(n+1)})$ is the archimedean set $(\Z,\theta^k)$, with $\theta^k:\Z\to\Z$, $\theta^k(x) = x + k(n+1)$.
\begin{lem}\label{lembB1}
$(i)$~The inclusion of categories $\dop\subset \Lambda$ induces an inclusion of sets $\Sigma_n^k\subset \Gamma_n^k$ and one has 
\begin{equation}\label{leftc}
\Gamma_n^k=\Aut_\Lambda(\Psi_k(\underline{n}))\Sigma_n^k.
\end{equation}
$(ii)$~The canonical decomposition $\Lambda=\Delta \mathbf C$ determines an inclusion of sets $\Delta_n^k\subset \Gamma_n^k$ and one has
\begin{equation}
\Gamma_n^k= \Delta_n^k\Aut_\Lambda(\underline n)
\end{equation}
$(iii)$~For any $\alpha \in \Sigma_{n+1}^k$ there exists a unique element $\alpha'\in \Delta_n^k$ such that $\alpha B_0=\sd_k(B_0)\alpha'$. Moreover the map $\alpha\mapsto \alpha'$ defines a  bijection of sets $\Sigma_{n+1}^k\simeq \Delta_n^k$.
\end{lem}
\proof $(i)$~Let $\alpha \in \Hom_\dop(n^*,\sdd_k(n^*))$, then it is easy to see that its image into the cyclic category $\tilde \alpha \in \Hom_\Lambda(\underline{n},\Psi_k(\underline{n}))$ belongs to $\Gamma_n^k$ if and only if $\alpha\in \Sigma_n^k$. In that case, the associated permutation of $\{0,\ldots, n\}$ is obtained by extending $\pee(\alpha)$ to a permutation of $\{0,1,\ldots, n\}$ by fixing $0$. One has $\Lambda=\mathbf C\dop$ and every morphism $\phi$ in $\Lambda$ uniquely decomposes as $\phi=t^a\circ \delta$, where $\delta$ is in $\dop$. For $\phi\in \Gamma_n^k$, one has $\phi=t^a\circ \delta$ where $\delta\in \Hom_\dop(n^*,\sdd_k(n^*))$, thus one gets \eqref{leftc}.

$(ii)$~By construction one has $\Gamma_n^k\cap \Delta=\Delta_n^k$. For $\alpha \in \Hom_\Lambda(\underline{n},\Psi_k(\underline{n}))$, the right multiplication $\alpha\mapsto \alpha \tau$ by an element  $\tau\in\Aut_\Lambda(n)$ does not affect the condition:  ``$ (\id_n^k\circ\alpha)_* \, \text{is bijective}$''  thus the result follows from the decomposition $\Lambda=\Delta \mathbf C$.

$(iii)$~For $\alpha \in \Hom_\dop((n+1)^*,(k(n+2)-1)^*)$  we let $\tilde \alpha$ its canonical lift to a non-decreasing map
$$
\tilde \alpha:\Z\to \Z, \  \ \tilde \alpha(x)=\alpha(x)\qqq x\in \{0, \ldots , n+1\}, \  \  \tilde \alpha(x+n+2)= \tilde \alpha(x)+k(n+2)\qqq x\in \Z.
$$
By hypothesis the map $\pee(\alpha)$ which to $x\in  \{1, \ldots , n+1\}$ associates $\alpha(x)$ modulo $n+2$ is a permutation of $ \{1, \ldots , n+1\}$. Thus the map $\phi:\Z\to \Z$, $\phi=\tilde \alpha\circ B_0$ is increasing and fulfills the following properties 
\begin{itemize}
\item $\phi(x+a(n+1))=\phi(x)+ka(n+2)$, $\forall x, a\in \Z$
\item $\phi(x)\neq 0$ modulo $n+2$, $\forall x\in \Z$.
\end{itemize}
We show that these conditions imply that there exists $\alpha'\in \Hom_\Delta([n],[k(n+1)-1])$ such that $\phi =\sd_k(B_0)\alpha'$, \ie the following diagram commutes
$$
\xymatrix@C=45pt@R=85pt{\{0,\ldots, n, n+1\} \  \  \ar[r]^(.45)\alpha &   \{0,\ldots,k (n+2)-1\}  \\ \{0,\ldots, n\}\  \ar[u]^{B_0}  \ar[r]^(.4){\alpha'}&\   \ar[u]^{\sd_k(B_0)} \{0,\ldots,k (n+1)-1\}}
$$
 One has $B_0\in \Hom_\Delta([n],[n+1])\subset \Hom_\Lambda(\underline n, \underline{ n+1})$ and the morphism $\sd_k(B_0)$ which belongs to $\Hom_\Delta(\sd_k([n]),\sd_k([n+1]))$ is also given by  \eqref{defnBo}. Thus the range of $\phi:\Z\to \Z$ is contained in the range of the injection $\sd_k(B_0)$ and for each $x\in \Z$ there exists a unique $y$ such that $\phi(x)= \sd_k(B_0)(y)$. Set $ \alpha'(x)=y$, then $\alpha':\Z\to \Z$ is a non-decreasing map which satisfies
\begin{itemize}
\item $\alpha'(x+a(n+1))=\alpha'(x)+ka(n+1)$, $\forall x, a\in \Z$
\item $\alpha'(x)\in \{0, \ldots , k(n+1)-1\}$, $\forall x\in \{0, \ldots , n\}$.
\end{itemize}
Thus $\alpha'\in \Hom_\Delta([n],[k(n+1)-1])$ and the map $(\ps_n^k\circ\alpha')_*$ is the same as the permutation $\pee(\alpha)$ of the set $ \{1, \ldots , n+1\}$ reindexed as $ \{0, \ldots , n\}$, thus $\alpha'\in \Delta_n^k$. 
Moreover the map $\alpha\mapsto \alpha'$ determines a  bijection $\Sigma_{n+1}^k\simeq\Delta_n^k$ since it is injective by construction and the cardinality of $\Sigma_{n+1}^k$, which is $k^{n+1}$ by Lemma \ref{lemparam}, is the same as the cardinality of $\Delta_n^k$ 
(which is by $(ii)$ equal to $1/n\times \#(\Gamma_n^k)=k^{n+1}$ using $(i)$). \qed 

\begin{thm} \label{thmBlamcom}
For any $n,k\geq 1$ one has  in $\Z[\lbt]$
\begin{equation}\label{Blamcom}
\Lambda_{n+1}^k B(n) =kB(n)  \Lambda_n^k
\end{equation}
\end{thm}
\proof
The left hand side of \eqref{Blamcom} is a sum of $(n+1)k^{n+1}$ terms which, up to sign, are of the form
$$
\ps_{n+1}^k\alpha B_0 \tau, \ \alpha \in \Sigma_{n+1}^k, ~\tau \in \Aut_\Lambda(\underline n).
$$
By Lemma \ref{lembB1} one has $\alpha B_0=\sd_k(B_0)\alpha'$, so that these terms can be rewritten as follows
\begin{equation}\label{lhs}
\ps_{n+1}^k\sd_k(B_0)\alpha' \tau, \ \alpha' \in \Delta_{n}^k, \tau \in \Aut_\Lambda(\underline n).
\end{equation}
The coefficient of $k$ in the right hand side of \eqref{Blamcom} is a sum of $(n+1)k^{n}$ terms which, up to sign, are of the form
$
B_0\,\xi\, \ps_n^k\, \beta, \ \beta \in \Sigma_{n}^k, \ \xi \in \Aut_\Lambda(\underline n)
$.

One also has $\xi \ps_n^k =\ps_n^k\tilde \xi$ where $\tilde \xi\in \Aut_\Lambda(\Psi_k(\underline{n}))$ is a lift of the cyclic permutation $\xi\in \Aut_\Lambda(n)$. Moreover $B_0 \ps_n^k =\ps_{n+1}^k\sd_k(B_0)$. Thus we conclude that 
$
B_0\,\xi\, \ps_n^k\, \beta=\ps_{n+1}^k\sd_k(B_0)\tilde \xi\,\beta
$. 
One  has 
$\tilde \xi\, \beta\in \Aut_\Lambda(\Psi_k(\underline{n}))\Sigma_n^k$ and using the coefficient $k$ in front of the right hand side of  \eqref{Blamcom} together with the equality $\ps_n^k \eta=\ps_n^k \eta'$ for any two of the $k$ lifts of $\xi$ one can rewrite the right hand side of  \eqref{Blamcom} as a sum of $(n+1)k^{n+1}$ terms which up to sign are of the form
\begin{equation}\label{rhs}
\ps_{n+1}^k\sd_k(B_0)\tilde \xi\,\beta, \  \ \tilde \xi\in \Aut_\Lambda(\Psi_k(\underline{n})), \,\beta\in \Sigma_n^k.
\end{equation}
Thus the validity of \eqref{Blamcom} follows from Lemma \ref{lembB1} which gives  $\Aut_\Lambda(\Psi_k(\underline{n}))\Sigma_n^k= \Gamma_n^k=\Delta_n^k\Aut_\Lambda(\underline n)$ and the fact that the signs in front of the terms in \eqref{lhs} and \eqref{rhs} are always given by the signature of the associated permutations.\qed

%B4
\subsection{$\lambda$-operations on the cyclic homology of epicyclic modules}

Let $E$ be an epicyclic module. Then its restriction to $\Lambda\subset \lbt$  is a cyclic module and one has the associated abelian groups \eqref{chains}
$$
C_n(E):=E(n)/V(n), \ \ V(n):=\oplus \,{\rm Im}(s_j).    
$$
The operations $\Lambda_n^k$ are meaningful on the quotient in view of the following

\begin{lem}\label{lamnkquot}
$
\Lambda_n^k\left(\oplus \,{\rm Im}(s_j)    \right)\subset \left(\oplus \,{\rm Im}(s_\ell)    \right), ~\forall j,\ell.
$
\end{lem}
\proof
For any $\sigma\in \Sigma_n^k$ and any degeneracy $s_j$, there exists $\alpha$ such that with $i=\pee(\sigma)(j+1)$, one has
$
\sigma\circ s_j=\sdd_k(s_{i-1})\circ \alpha
$.
The conclusion follows using $\ps_n^k \,\sdd_k(s_{i-1})=s_{i-1}\ps_n^k$.\qed

\begin{lem} \label{lemlambdaop} $(i)$~The following equality defines, for any integer $k\ge 1$, an endomorphism of the $(b,B)$-bicomplex of $E$:
\begin{equation}\label{theta}
\theta(k)\xi=k^{-\alpha}\Lambda_{\alpha-\beta}^k\,\xi\qqq \xi \in C^{\alpha,\beta}
\end{equation}
$(ii)$~The translation $(\alpha,\beta)\mapsto (\alpha+1,\beta+1)$ defines an endomorphism $S$ 
of the $(b,B)$-bicomplex and one has 
$
S\theta(k)=k\theta(k)S.
$
\end{lem}
\proof $(i)$~The commutation of $\theta(k)$ with $b$ and $B$ follows from Proposition \ref{propblamcom} and Theorem \ref{thmBlamcom} respectively. More precisely one has for $\xi \in C^{\alpha,\beta}$, $B\xi \in C^{\alpha+1,\beta}$ and
$$
B\theta(k)\xi=Bk^{-\alpha}\Lambda_{\alpha-\beta}^k\,\xi=k^{-\alpha-1}kB\Lambda_{\alpha-\beta}^k\,\xi=
k^{-(\alpha+1)}\Lambda_{\alpha+1-\beta}^k\,B\xi=\theta(k)B\xi.
$$
$(ii)$~Follows directly from the definition of $\theta(k)$ as in \eqref{theta}.\qed

The endomorphisms $\theta(k)$ as in \eqref{theta} define, at the categorical level, the $\lambda$-operations in homology. Lemma \ref{lemlambdaop}  immediately implies the following well-known result (\cf \cite{Loday}, Ex. 6.4.5) 
\begin{thm}\label{thmlambdaop}
Let $E$ be an epicyclic module. The $\lambda$-operations define an action of $\N^\times$ on $HC_n(E)$
given by the induced action of the operators $\theta(k)$. 
\end{thm}

When the epicyclic module $E$ factors through  $\fin$ 
one has, after tensorisation of  the abelian groups by $\Q$, a  decomposition (\cf \cite{Loday} Thm. 6.4.5)  
$$
    HC_n(E)=\bigoplus_{j\ge 0} HC_n^{(j)}(E)
$$
which diagonalizes the action of $\nt$ as a sum of the representations given by the characters $\nt\ni k\mapsto k^j$. This decomposition of an epicyclic module does not hold in general. Let $E$ be an epicyclic module and $\rho:\nt\to V$ be a representation of $\N^\times$, then the tensor product $E\otimes V$  endowed with the maps
$$
(E\otimes V)(f)=E(f)\otimes \rho(\mmod(f)):E(n)\otimes V\to E(m)\otimes V\qqq f\in \Hom_\lbt(\underline n,\underline m)
$$ 
is still an epicyclic module. When working with vector spaces over a field, this construction tensors the cyclic homology by $V$ and replaces the action of $\nt$ on $HC_n(E)$ by its tensor product with $\rho$. Thus this twisting process generates epicyclic modules with arbitrary weights.

\end{document}